\documentclass[a4paper,11pt]{article}
\usepackage{stmaryrd}
\usepackage{amsfonts}
\usepackage{bbm}
\usepackage{float}
\usepackage{amscd}
\usepackage{mathrsfs}
\usepackage{latexsym,amssymb,amsmath,amscd,amscd,amsthm,amsxtra,xypic}
\usepackage[dvips]{graphicx}
\usepackage[utf8]{inputenc}
\usepackage[T1]{fontenc}
\usepackage{lmodern}
\usepackage{amssymb}
\usepackage[all]{xy}
\usepackage{nicefrac,mathtools,enumitem}
\usepackage{microtype}
\usepackage{lipsum}
\usepackage{mathtools}
\usepackage{xfrac}

\usepackage{geometry}
\usepackage{hyperref}

\usepackage{times}
\usepackage{amsthm}
\usepackage{amsmath}
\usepackage{amsfonts}
\usepackage{amssymb}
\usepackage{mathdots}
\usepackage{color}
\usepackage{mathrsfs}
\usepackage{stmaryrd}
\SetSymbolFont{stmry}{bold}{U}{stmry}{m}{n}

\usepackage{bm}
\usepackage{tikz-cd}
\usepackage{tikz}
\usetikzlibrary{matrix,arrows,decorations.pathmorphing}
\usepackage[all]{xy}
\usepackage{graphicx}
\usepackage{subfigure}
\usepackage[toc,page]{appendix}
\usepackage[usestackEOL]{stackengine}

\newcommand{\I}{{\mathrm{i}}}

\usepackage{theoremref}
\newtheorem{thm}{Theorem}[section]
\newtheorem{lem}[thm]{Lemma}
\newtheorem{cor}[thm]{Corollary}
\newtheorem{pro}[thm]{Proposition}

\theoremstyle{definition}

\newtheorem{rmk}[thm]{Remark}
\newtheorem{defi}[thm]{Definition}

\newcommand {\Sect}{{\rm Sect}}

\newcommand{\be }{\begin{equation}}
\newcommand{\ee }{\end{equation}}

\newcommand{\pf}{\noindent{\bf Proof.}\ }

\newcommand{\h}{\mathfrak h}

\def\qed{\hfill ~\vrule height6pt width6pt depth0pt}

\newcommand{\br}[1]{   [ \cdot,    \cdot  ]   }

\newcommand{\g}{\mathfrak g}

\newcommand{\gl}{\mathfrak {gl}}

\title{WKB approximation, crystals and combinatorics of Young tableaux}
\author{Xiaomeng Xu}
\date{}

\newcommand{\Addresses}{{
  \bigskip
  \footnotesize
\noindent \textsc{School of Mathematical Sciences \& Beijing International Center
for Mathematical Research, Peking University, Beijing 100871, China}\par\nopagebreak
  \textit{E-mail address}: \texttt{xxu@bicmr.pku.edu.cn}
}}
\begin{document}

\maketitle
\begin{abstract}
In this paper, we show how various combinatorial algorithms of Young tableaux naturally arise from the WKB approximation of the connection matrix of quantum confluent hypergeometric equation, including the Schensted bumping operator, the Robinson-Schensted correspondence, and the Littlewood-Richardson rule.
\end{abstract}

\section{Introduction}

Young tableaux play an essential role in the representation theory of the general linear Lie algebra. Regarding the combinatorics of Young tableaux, there are many important algorithms that are closely related to representation theory, for example, the Schensted bumping rule, the
Sch\"utzenberger involution, the Robinson-Schensted correspondence and the Littlewood-Richardson rule (see, e.g., \cite{Fu, Sun}). These algorithms have been understood using the crystal structures of the quantum group
$U_q(\frak{gl}_n)$, see, e.g., \cite{Kwon} for such an exposition.

In a previous paper, we proved that the WKB approximation of the Stokes matrices of the quantum
confluent hypergeometric equation gives rise to $\frak{gl}_n$-crystals realized on the semistandard Young tableaux. The purpose of this paper is to deepen the relation between the WKB approximation and representation theory, by realizing various combinatorics of Young tableaux via the WKB approximation as shown in the following table.

\begin{table}[h!]
\caption{WKB approximation in the Stokes phenomenon of equation \eqref{introqeq} and combinatorics of Young tableaux}
\begin{tabular}{ |c|c|c| } 
 \hline
1 & WKB approximation of Stokes matrices & ${\gl}_n$-crystals on semistandard Young Tableaux \\  \hline 
2 & WKB approximation of connection matrices & Robinson-Schensted algorithm
\\
 \hline
3 & WKB approximation of coproduct & { Littlewood–Richardson rule }
\\ \hline
4 & Wall-crossing in the WKB approximation &  Sch\"utzenberger involution
\\ \hline
%5 & WKB of degenerate $u={{\rm diag}(\underbrace{0,...,0}_{m}, {u_{m+1},...,u_n})}$ & Branching rules for $\frak{gl}_m\subset \frak{gl}_n$ 
%\\ \hline
\end{tabular}
\end{table}

Let us take the Lie algebra ${\frak {gl}_n}$ over the field of complex numbers, and its universal enveloping algebra $U({\frak {gl}}_n)$ generated by $\{e_{ij}\}_{1\le i,j\le n}$ subject to the relation $[e_{ij},e_{kl}]=\delta_{jk}e_{il}-\delta_{li} e_{kj}$. Let us take the $n\times n$ matrix $\mathcal{T}=(\mathcal{T}_{ij})$ with entries valued in $U({\frak {gl}}_n)$
\begin{eqnarray*}
\mathcal{T}_{ij}=e_{ji}, \ \ \ \ \ \text{for} \ 1\le i,j\le n.
\end{eqnarray*}
Let $\h_{\rm reg}(\mathbb{R})$ denote the set of $n\times n$ diagonal matrices with distinct eigenvalues. Given any finite-dimensional representation $W$ of ${\gl}_n$, let us consider the quantum confluent hypergeometric system
\begin{eqnarray}\label{introqeq}
\frac{dF}{dz}=h\Big( \I u+\frac{1}{2\pi\I}\frac{\mathcal{T}}{z}\Big)\cdot F,
\end{eqnarray}
for $F(z)\in {\rm End}(W)\otimes {\rm End}(\mathbb{C}^n)$ an $n\times n$ matrix function with entries in ${\rm End}(W)$. Here $\I=\sqrt{-1}$, $h>0$ is a positive real parameter, $u={\rm diag}(u_1,...,u_n)\in\h_{\rm reg}(\mathbb{R})$ is seen as a $n\times n$ matrix with scalar entries in $U({\frak {gl}}_n)$, and the action of the coefficient matrix on $F(z)$ is given by matrix multiplication and the representation of ${\frak {gl}}_n$.

The equation \eqref{introqeq} has a unique formal solution $\hat{F}$ around
$z = \infty$. The standard theory of resummation states that there exist certain sectorial regions around $z=\infty$, such that on each of these sectors there is a unique (therefore canonical) holomorphic solution with the prescribed asymptotics $\hat{F}$. These solutions are in general different (that reflects the Stokes phenomenon), and the transition between them can be measured by a pair of Stokes matrices $S_\pm(u)^W\in{\rm End}(W)\otimes{\rm End}(\mathbb{C}^n)$. In the meanwhile, the equation \eqref{introqeq} has a canonical holomorphic solution around
$z = 0$, and the comparison of solutions at $0$ and $\infty$ defines a connection matrix $C(u)^W\in {\rm End}(W)\otimes{\rm End}(\mathbb{C}^n)$. See Section \ref{solmono} for more details. 
It was proved in \cite{Xu2} that
\begin{thm}\cite{Xu2}\label{quantumgpucat}
For any fixed positive real number $h$, the map (with $q=e^{h/2}$)
\begin{equation*}
\begin{split}
\mathcal{S}_q(u): U_q(\frak{gl}_n)&\rightarrow {\rm End}(W)~;\\
e_i&\mapsto \frac{ 1}{q-q^{-1}} q^{-e_{i+1,i+1}}\cdot S_-(u)^W_{i+1,i}\cdot q^{e_{ii}}, \\
f_i&\mapsto \frac{1}{q^{-1}-q} S_+(u)^W_{i,i+1}, \\
q^{h_i}&\mapsto q^{-e_{ii}}
\end{split}
\end{equation*}
defines a representation of the Drinfeld-Jimbo quantum group $U_q(\frak{gl}_n)$ on the vector space $W$. Here $U_q(\frak{gl}_n)$ is a unital associative algebra with generators $q^{\pm h_i}, e_j, f_j,$ $1\le j\le n-1, 1\le i\le n$, see Definition \ref{defqg}.
\end{thm}
This theorem motivates us to further explore the relation between the Stokes phenomenon and representation theory. However, compared to the representation theory of quantum groups, less is known about the quantitative analysis of the Stokes matrices. Thus, to deepen the relation, we first need to develop some analysis tool. 
In general, the Stokes and connection matrices are (new) highly transcendental matrix functions depending on the chosen representation $W$ and the variables $u={\rm diag}(u_1,...,u_n)$. For fixed $W$ they are real analytic functions of $u\in \h_{\rm reg}(\mathbb{R})$ and have singularities along the hyperplanes
\[\h(\mathbb{R})\setminus \h_{\rm reg}(\mathbb{R})=\{(u_1,...,u_n)\in \mathbb{R}^n~|~u_i = u_j, \text{for some } i\ne j \}.\]
Therefore, if we want to study the dependence of the Stokes matrices on the representation spaces, we can take a special $u$ to minimize analytical complexity as much as possible. Our strategy is to take the reference point $u$ as one of the critical points. 

From the viewpoint of analysis, to understand the transcendental functions, one should study their expansions around critical points. This problem was studied in our previous work \cite{Xu}: the leading terms in the expansion or equivalently the regularized limits of the Stokes matrices and connection matrices were given, as the components $u_i$ of $u\in\h_{\rm reg}(\mathbb{R})$ collapse at a comparable speed. The prescription of the regularized limits is controlled by the geometry of the De Concini-Procesi space $\widetilde{\frak t_{\rm reg}}(\mathbb{R})$. Here, roughly speaking, the De Concini-Procesi space replaces the set $\h(\mathbb{R})\setminus \h_{\rm reg}(\mathbb{R})$ of hyperplanes in $\h(\mathbb{R})$ by a divisor with normal crossings and leaves the complement $\h_{\rm reg}(\mathbb{R})$ of these subspaces
unchanged. In particular, denote by $u_{\rm cat}$ the limit
 \begin{equation}\label{limit}
    u_2-u_1\rightarrow 0+ \text{ and } \frac{u_j-u_{j-1}}{u_{j-1}-u_{j-2}}\rightarrow +\infty \text{ for all } j=3,...,n,
\end{equation}
(that is understood as a point in the $0$-dimensional stratum of $\widetilde{\frak h_{\rm reg}}(\mathbb{R})$, called a caterpillar \cite{Sp}), the regularized limits of the Stokes and connection matrices as $u\rightarrow u_{\rm cat}$ are denoted by $S_\pm(u_{\rm cat})^W$ and $C(u_{\rm cat})^W\in{\rm End}(W)\otimes{\rm End}(\mathbb{C}^n)$ respectively. They can be understood as the monodromy matrices of the equation at the limit point $u_{\rm cat}$. The advantage of these matrices at $u_{\rm cat}$ is that they have closed formula! See Sections \ref{sec:qStokescat} and \ref{proofofqS} for their definition and explicit expression, respectively.

Now, from the viewpoint of algebra, for all $u$ (therefore for the limit $u_{\rm cat}$), the Stokes matrices $S_\pm(u)^W$ give rise to representations of quantum groups as in Theorem \ref{quantumgpucat}. The Stokes matrices $S_\pm(u_{\rm cat})^W$ at the caterpillar point $u_{\rm cat}$ already capture the algebraic properties shared by the Stokes matrices at a generic point $u$. Since $S_\pm(u_{\rm cat})^W$ and $C(u_{\rm cat})^W$ have closed formula, one can then simply use them as an analytic model to unveil various correspondences between the transcendental structures in Stokes phenomenon (at a second order pole) and combinatorial structures in representation theory. In the following, we show such examples.

\subsection{${\gl}_n$-crystals arising from the WKB approximation of Stokes matrices}

Given any finite dimensional irreducible representation $L(\lambda)$ of ${\gl}_n$ with a highest weight $\lambda$, we take a set $E_{GT}(\lambda)$ of Gelfand-Tsetlin (GT) orthonormal basis $\{\xi_\Lambda\}\in L(\lambda)$, parameterized by the Gelfand-Tsetlin patterns $\Lambda \in GT_n(\lambda)$. See Section \ref{GZbasis}.

\begin{thm}\cite{Xu}\label{WKBthm}
For each $k=1,...,n-1$, there exists canonical operators $\widetilde{e_k}$ and $\widetilde{f_k}$ acting on the finite set $E_{GT}(\lambda)$, and real valued functions $c_{ki}(\Lambda)$, $\theta_{ki}(\Lambda)$ with $i=1,2$ such that for any generic element $\xi_\Lambda\in E_{GT}(\lambda)$,
\begin{eqnarray*}
\mathop{\rm lim}\limits_{h\rightarrow +\infty}\left( S_{+}(u_{\rm cat})^{L(\lambda)}_{k,k+1}\cdot e^{h c_{k1}(\Lambda) +\I  h\theta_{k1}(\Lambda)} \xi_\Lambda\right)=\widetilde{f_k}(\xi_\Lambda),\\
\mathop{\rm lim}\limits_{h\rightarrow +\infty}\left( S_{-}(u_{\rm cat})^{L(\lambda)}_{k+1,k}\cdot e^{hc_{k2}(\Lambda) +\I h  \theta_{k2}(\Lambda)} \xi_\Lambda\right)=\widetilde{e_k}(\xi_\Lambda).
\end{eqnarray*}
(We refer Section \ref{sec:WKB} for the notion of generic basis elements). Furthermore, the set $E_{GT}(\lambda)$ equipped with the operators $\widetilde{e_k}$ and $\widetilde{f_k}$ is a ${\gl}_n$-crystal.
\end{thm}

Under the bijection between the set $B_n(\lambda)$ of semistandard Young tableaux and the set $GT_n(\lambda)$ of Geland-Testlin patterns (see Section \ref{GTY}), the $\frak{gl}_n$-crystal induced from the WKB approximation coincides with the Young tableaux realization of $\frak{gl}_n$-crystals \cite{KN}. We note that the WKB approximation of the Stokes matrices of the (classical) confluent hypergeometric equations has been studied in \cite{ANXZ} with a relation to cluster algebras.

\subsection{Schensted deletion operator arising from the WKB approximation of connection matrices}

Let $V=\mathbb{C}^n$ be the $n$ dimensional irreducible representation of $\frak{gl}_n$, and
let $\{v_i\}_{i=1,...,n}$ be its natural basis. 
Let us take a partition $\lambda=(\lambda_1,...,\lambda_n)$, seen as a Young diagram.
The decomposition of the tensor product $L(\lambda)\otimes V$ of the representation $L(\lambda)$ and the natural representation $V$ is
\begin{equation}
    L(\lambda)\otimes \mathbb{C}^n=\underset{\mu}{\bigoplus} L(\mu)
\end{equation}
where $\mu$ runs over Young diagrams obtained by adding one node to $\lambda$. We denote by $\lambda\stackrel{a}{\longleftarrow} \mu$, if the node is added at the $a$-th row. That is $\mu=(\mu_1,...,\mu_n)$ with $\mu_a=\lambda_a+1$, and $\mu_i=\lambda_i$ for $i\ne a$. 

Recall that the connection matrix $C(u_{\rm cat})^{L(\lambda)}\in {\rm End}(L(\lambda))\otimes {\rm End}(\mathbb{C}^n)$. In Section \ref{SfromWKB}, we prove the main theorem of this paper, i.e., 
the following transcendental realization of the Schensted deletion operation.
\begin{thm}\label{WKBC0}
Given any basis vector $\xi_{\Lambda}\in L(\mu)\subset  L(\lambda)\otimes \mathbb{C}^n$ with $\lambda\stackrel{a}{\longleftarrow} \mu$,
there exists real number $\theta(\Lambda)$ such that 
\begin{eqnarray*}
\mathop{\rm lim}\limits_{h\rightarrow +\infty} C(u_{\rm cat})^{L(\lambda)}\cdot e^{\I h \theta(\Lambda)} \xi_\Lambda =\xi_{\leftarrow_a\Lambda}\otimes v_j.
\end{eqnarray*}
Here $\leftarrow_a\Lambda$ is the Gelfand-Testlin pattern obtained from $\Lambda$ by the Schensted deletion operation (see Definition \ref{delop}), and $\xi_{\leftarrow_a\Lambda}\in L(\lambda)$ is the corresponding basis vector.
\end{thm}
For a constant vector $w\in L(\lambda)$ (that is independent of $h$) and a $h$ dependent action $\rho(h)\in {\rm End}(L(\lambda))$, we will write 
\[\mathcal{WKB}(\rho(h)\cdot v)= v'\]
if there exist some complex number $c$ and a nonzero constant vector $v'$ such that $\mathop{\rm lim}\limits_{h\rightarrow +\infty} \rho(h)\cdot e^{h c} v=v'$. For example, Theorem \ref{WKBthm} and Theorem \ref{WKBC0} state that 
\[\mathcal{WKB}\left( S_{+}(u_{\rm cat})^{L(\lambda)}_{k,k+1}\cdot \xi_\Lambda\right)=\widetilde{f_k}(\xi_\Lambda) \text{  and  } \mathcal{WKB}\left(C(u_{\rm cat})^{L(\lambda)}\cdot \xi_\Lambda\right) =\xi_{\leftarrow_a\Lambda}\otimes v_j. \]
\subsection{Robinson–Schensted correspondence}
Recall that $V$ is the $n$ dimensional irreducible representation. Denote by ${\rm Id}\in {\rm End}(V)$ the identity and then 
\begin{equation}\label{mC}
    {\rm Id}^{N-k-1}\otimes C(u_{\rm cat})^{V^{\otimes k}}\in {\rm End}(V^{\otimes N}) \ \text{   for } k=1,...,N-1.
\end{equation} The tensor product has an irreducible decomposition 
$V^{\otimes N}=\underset{W}{\bigoplus} V_W$ where $W$ is a standard tableau with $N$ nodes. 
We consider the action of the operators in \eqref{mC} on the GT basis of the irreducible representation $V_W\subset V^{\otimes N}$. By repeatedly utilizing Theorem \ref{WKBC0}, we get
\begin{thm}\label{WKBmC0}
For any Gelfand-Tsetlin basis vector $\xi_\Lambda\in V_W\subset V^{\otimes N}$, there exists a real function $\theta_N(\Lambda)$ such that
\begin{eqnarray*}
\mathop{\rm lim}\limits_{h\rightarrow +\infty}\prod_{k=1}^{N-1} \left({\rm Id}^{N-k-1}\otimes C(u_{\rm cat})^{V^{\otimes k}}\right)\cdot q^{\I h\theta_N(\Lambda)}\xi_\Lambda=v_{w_1}\otimes \cdots \otimes v_{w_N},
\end{eqnarray*}
where the correspondence $\xi_\Lambda\rightarrow v_{w_1}\otimes \cdots \otimes v_{w_N}$ is given by the Robinson-Shensted 
correspondence.
\end{thm}
Here the Robinson-Shensted 
correspondence is a bijection between certain sequences of numbers and pairs of Young tableaux, defined by means of a combinatorial algorithm, see Section \ref{RScorr} for more details. There are algebraic and geometric interpretations of the Robinson-Schensted algorithm, for example, via the primitive ideals in enveloping algebras (see \cite[Theorem 6.5]{Vog}), the cells
in Coxeter groups \cite{KL}, subvarieties of the flag manifolds \cite{St}, and the representation theory of the quantum group $U_q(\frak{gl}_n)$ at $q = 0$ \cite{DJM}. As far as we know, the above theorems give the first transcendental realization.

\subsection{Littlewood–Richardson rule and the WKB approximation of coproduct}
Let us take the irreducible decomposition of the tensor product of two irreducible representations $L(\mu)$ and $L(\nu)$
\begin{equation}
L(\mu)\otimes L(\nu) =\bigoplus_\lambda m^\lambda_{\mu\nu} L(\lambda),
\end{equation}
where $m^\lambda_{\mu\nu}$ are the Littlewood–Richardson coefficients.

There are two different actions on the tensor product. On the one hand, 
we can consider the entries $S_{+}(u_{\rm cat})^{L(\mu)\otimes L(\nu)}_{k,k+1}$ and $S_{+}(u_{\rm cat})^{L(\mu)\otimes L(\nu)}_{k,k+1}$ of the Stokes matrices of the linear system \eqref{introqeq} associated to the tensor representation $W=L(\mu)\otimes L(\nu)$. Their actions on $L(\mu)\otimes L(\nu)$ can be restricted to each component $L(\lambda)\subset L(\mu)\otimes L(\nu)$.
On the other hand, we can consider the actions of
\begin{align}\label{copro1}
\Delta(S_{+}(u_{\rm cat})_{k,k+1})^{L(\mu)\otimes L(\nu)}:&=q^{e_{kk}-e_{k+1,k+1}}\otimes S_{+}(u_{\rm cat})^{L(\nu)}_{k,k+1}+S_{+}(u_{\rm cat})^{L(\mu)}_{k,k+1}\otimes 1\\  \label{copro2}
\Delta(S_{-}(u_{\rm cat})_{k+1,k})^{L(\mu)\otimes L(\nu)}:&=q^{e_{k+1,k+1}-e_{kk}} \otimes S_{-}(u_{\rm cat})^{L(\nu)}_{k+1,k}+S_{-}(u_{\rm cat})^{L(\mu)}_{k+1,k}\otimes q^{2e_{k+1,k+1}-2e_{kk}},
\end{align} 
on the tensor product $L(\mu)\otimes L(\nu) $. The operators \eqref{copro1}-\eqref{copro2} are induced from the coproduct of the quantum group $U_q(\frak{gl}_n)$, see Section \ref{qgStokes}.

The difference of the two actions in the WKB approximation can be described in a combinatorial way:

\begin{thm}
Let $\xi_{\Lambda_1}$ and $\xi_{\Lambda_2}$ be any generic GT basis vectors of $L(\mu)$ and $L(\nu)$ respectively. Then we have
\begin{align*}
&LR\left(\mathcal{WKB}\left( \Delta(S_{+}(u_{\rm cat}))^{L(\mu)\otimes L(\nu) }_{k,k+1}\cdot (\xi_{\Lambda_1}\otimes \xi_{\Lambda_2})\right)\right)=\mathcal{WKB}\left( (S_{+}(u_{\rm cat})^{L(\mu)\otimes L(\nu) })_{k,k+1}\cdot LR(\xi_{\Lambda_1}\otimes \xi_{\Lambda_2}) \right),\\
&LR\left(\mathcal{WKB}\left( \Delta(S_{-}(u_{\rm cat}))^{L(\mu)\otimes L(\nu) }_{k+1,k}\cdot (\xi_{\Lambda_1}\otimes \xi_{\Lambda_2})\right)\right)=\mathcal{WKB}\left( (S_{-}(u_{\rm cat})^{L(\mu)\otimes L(\nu) })_{k+1,k}\cdot LR(\xi_{\Lambda_1}\otimes \xi_{\Lambda_2}) \right),
\end{align*}
where the map $LR:L(\mu)\otimes L(\nu)\rightarrow \bigoplus_\lambda m^\lambda_{\mu\nu} L(\lambda)$ is given by the Littlewood–Richardson rule (see Section \ref{LRWKB} for more details).
\end{thm}
This theorem can be used to compute the WKB approximation of the Stokes matrices $S_\pm(u_{\rm cat})$ with respect to the tensor product of multiple irreducible representations.

\subsection{Sch\"utzenberger involution and Wall-crossing}
Though we will not pursue the following viewpoint further here, it is worth recalling the wall-crossing of the WKB approximation of the Stokes matrices along codimension one boundaries in the compactified space $\widehat{\h_{\rm reg}}(\mathbb{R})$.

The connected components of $\h_{\rm reg}(\mathbb{R})$ are
$U_{\sigma}=\{u_{\sigma(1)}<\cdots <u_{\sigma(n)}\}$ for some permutation $\sigma\in S_n$. 
For each $1\le i\le n$, let $\tau_i\in S_n$ be the permutation reversing the segment $[1,...,i]$, and \[U_{\tau_i}=\{u_{i}<u_{i-1}<\cdots <u_1<u_{i+1}\cdots <u_{n}\}\] the corresponding connected component. If we think of $u_{\rm cat}$ as a point on $\widehat{\h_{\rm reg}}(\mathbb{R})$, then for all $i=1,...,n$, it lies at the closure of ${U_{\tau_i}}$ in the compactification $\widehat{\h_{\rm reg}}(\mathbb{R})$: $u_{\rm cat}$ also represents the following limit point of $u\in U_{\tau_i}$,
 \begin{equation}\label{limiti}
     u_{\tau_i(2)}-u_{\tau_i(1)}\rightarrow 0+ \text{ and } \frac{u_{\tau_i(j)}-u_{\tau_i(j-1)}}{u_{\tau_i(j-1)}-u_{\tau_i(j-2)}}\rightarrow +\infty \text{ for all } j=3,...,n.
\end{equation}

Accordingly, the regularized limit of $S_{\pm}(u)^{L(\lambda)}$, as $u\rightarrow u_{\rm cat}$ from the connected component $U_{\tau_i}$, are denoted by $S_{\pm}^{\tau_i}(u_{\rm cat})^{L(\lambda)}\in{\rm End}(L(\lambda))\otimes {\rm End}(\mathbb{C}^n)$. Thus, $S_{\pm}(u_{\rm cat})$ are simply $S_{\pm}^{\rm id}(u_{\rm cat})$. Following \cite{Xu}, these regularized limits $S_{\pm}^{\tau_i}(u_{\rm cat})^{L(\lambda)}$, as $u\rightarrow u_{\rm cat}$ from different connected components, are different, and such phenomenon is called the wall-crossing (along the low dimensional stratum/boundary $\widetilde{\h_{\rm reg}}(\mathbb{R})\setminus {\h_{\rm reg}}(\mathbb{R})$ of the compactified space $\widetilde{\h_{\rm reg}}(\mathbb{R})$). In our previous paper \cite{Xu}, it was proved that the partial Sch\"utzenberger involutions, on the set of GT basis of $L(\lambda)$ or equivalently the set of semistandard Young tableaux, measure the difference of the WKB approximation of $S_{\pm}(u_{\rm cat})^{L(\lambda)}$ and $S_{\pm}^{\tau_i}(u_{\rm cat})^{L(\lambda)}$ for all $i=2,...,n$.

In the end, we stress that it is the closed formula/quantitative analysis of the Stokes matrices that makes it possible to discover the transcendental realizations of the various combinatorial structures.

\section{Crystals and Young tableaux}

\subsection{Crystals and tensor products}\label{sec:crystal}
Let $\g$ be a semisimple Lie algebra with a Cartan datum $(A,\Delta_+=\{\alpha_i\}_{i\in I},\Delta^\vee_+=\{\alpha_i^\vee\}_{i\in I}, P, P^\vee)$ be a Cartan datum, where $P\subset \h^*$ denotes the weight lattice, $I$ denotes the set of vertices of its Dynkin diagram, $\alpha_i\in I$ denote its simple
roots, and $\alpha^\vee_i$ the simple coroots.

\begin{defi}\cite{Ka3,Ka4}\label{WKBdata}
A $\g$-crystal is a finite set $B$ along with maps
\begin{eqnarray*}
wt&:& B\rightarrow P,\\
\tilde{e}_i,\tilde{f}_i&:& B\rightarrow B\cup \{0\}, \hspace{2mm} i\in I,\\
\varepsilon_i,\phi_i&:& B\rightarrow \mathbb{Z}\cup \{-\infty\}, \hspace{2mm} i\in I,
\end{eqnarray*}
satisfying for all $b,b'\in B$, and $i\in I$,
\begin{itemize}
    \item $\tilde{f}_i(b)=b'$ if and if $b=\tilde{e}_i(b')$, in which case
    \[wt(b')=wt(b)-\alpha_i, \hspace{3mm} \varepsilon_i(b')=\varepsilon_i(b)+1, \hspace{3mm} \phi_i(b')=\phi_i(b)-1.\]
    \item $\phi_i(b)=\varepsilon_i(b)+\langle wt(b),\alpha_i^\vee\rangle$, and if $\phi_i(b)=\varepsilon_i(b)=-\infty$, then $\tilde{e}_i(b)=\tilde{f}_i(b)=0.$
\end{itemize}
The map $wt$ is called the weight map, $\tilde{e}_i$ and $\tilde{f}_i$ are called Kashiwara operators or crystal operators.
\end{defi}

In this paper, we focus on the case $\g=\frak{gl}_n$.
For $n\ge 2$, let $P_n = \oplus _{i=1,..., n} \mathbb{Z} \epsilon_i$ be the free abelian group with a basis $\{\epsilon_1,..., \epsilon_n\}$,
 called the weight lattice. There is a natural symmetric bilinear form on $P_n$ given by $(\epsilon_i, \epsilon_j)=\delta_{ij}$. Let $P_n^\vee={\rm Hom}(P_n,\mathbb{Z})$ be the dual weight lattice, and let $\langle,\rangle$ denote the canonical pairing on $P_n^\vee \times P_n$. Set $I_n=\{1,...,n-1\}$. Then, for each $i\in I_n$ the simple root $\alpha_i$ is given by $\alpha_i=\epsilon_i-\epsilon_{i+1}$, and the simple coroot $h_i$ is defined to be the unique element in $P_n^\vee$ such that $\langle h_i,\lambda\rangle =(\alpha_i,\lambda)$ for all $\lambda\in P_n$. Note that the matrix $(\langle h_j,\alpha_i\rangle)_{i,j\in I_n}$
is the Cartan matrix of $A_{n-1}$ type.

\begin{defi}\label{crytensor}
Let $B_1$ and $B_2$ be two crystals. The tensor product $B_1\otimes B_2$ is the crystal with the underlying set $B_1\times B_2$ (the Cartesian product) and structure maps
 \begin{eqnarray*}
wt(b_1,b_2)&=&wt(b_1)+wt(b_2),\\
\tilde{e}_i(b_1,b_2)&=&\left\{
          \begin{array}{lr}
             (b_1,\tilde{e}_i(b_2)),   & \text{if} \ \varepsilon_i(b_1) \le \phi_i(b_2)  \\
           (\tilde{e}_i(b_1),b_2), & \text{otherwise}
             \end{array}\right. \\
\tilde{f}_i(b_1,b_2)&=&\left\{
          \begin{array}{lr}
             (b_1, \tilde{f}_i(b_2)),   & \text{if} \ \varepsilon_i(b_1)< \phi_i(b_2)  \\
           (\tilde{f}_i(b_1),b_2), & \text{otherwise}.
             \end{array}
\right. 
\end{eqnarray*}
\end{defi}

\subsection{Crystal of words}
For $n\ge 2$, let 
\[B_n=\{1<2<\cdots<n\}\]
be a linearly ordered set. Then $B_n$ is a $\frak{gl}_n$-crystal with the crystal operators
\begin{align*}
    &1\stackrel{\tilde{f_1}}{\longleftarrow} 2 \stackrel{\tilde{f_2}}{\longleftarrow}\cdots \stackrel{\tilde{f}_{n-1}}{\longleftarrow} n,\\
   & 1\stackrel{\tilde{e_1}}{\longrightarrow} 2 \stackrel{\tilde{e_2}}{\longrightarrow}\cdots \stackrel{\tilde{e}_{n-1}}{\longrightarrow} n, 
\end{align*}
and for $b\in B_n$, $wt(b)=\epsilon_b$. 

Let $\mathcal{W}_n$ be the set of all finite words with letters in $B_n$. That is
\[\mathcal{W}_n=\left(\oplus_{r\ge 1} B_n^{\otimes r}\right)\oplus \{\emptyset\},\]
 where $\emptyset$ denotes the empty word. Then by the tensor product rule of the crystals, $\mathcal{W}_n$ is a $\frak{gl}_n$-crystal,  with $wt(\emptyset) = 0$, $\tilde{e_i}\emptyset=\tilde{f_i} \emptyset=0$, and $\varepsilon_i(\emptyset) = \phi_i(\emptyset) = 0$ for all $i\in I_n$.

\subsection{Crystal of Young tableaux}
A partition of length $n$ is a nonincreasing sequence of nonnegative integers $\lambda = (\lambda_1,...,\lambda_n)$. A partition $\lambda$ is identified with a Young diagram, which is a collection of nodes (or boxes) in left-justified rows with $\lambda_k$ nodes in the $k$-th row numbered from top to bottom. We denote by $\mathcal{P}_n$ the set of
 all partitions of length $n$.

\begin{defi}
A Young tableau $YT$ is obtained by filling a Young diagram $\lambda$ with nonnegative integers. Then $\lambda$ is called the shape of $YT$, and $YT$ is called semistandard if the entries in each row are weakly increasing from left to right. If all the entries in $YT$ are distinct, then $YT$ is called standard. 
\end{defi}
For $\lambda\in \mathcal{P}_n$, let $B_n(\lambda)$ be the set of all semistandard Young tableaux of shape $\lambda$ with entries in $B_n =\{1,...,n\}$. 
For each semistandard tableau $YT$, let $\omega(YT)$ denote the word obtained by reading the entries of $YT$ column by column from right to left, and in each column from top to bottom. A semistandard tableau $YT$ can be uniquely recovered from the associated word $\omega(YT)$. In this way, we get a natural embedding of the set of semistandard tableau into the set of words \[\iota: B_n(\lambda)\rightarrow\mathcal{W}_n~;~YT\mapsto \omega(YT).\] 

The set $B_n(\lambda)$, seen as a subset of $\mathcal{W}_n$, inherits a $\frak{gl}_n$-crystal structure from $\mathcal{W}_n$: following \cite{KN}, $B_n(\lambda)$ together with $\emptyset$ is stable under the crystal operatos $\tilde{e}_i$ and $\tilde{f}_i$ of $\mathcal{W}_n$. And the weight of $YT\in B_n(\lambda)$ is given by $wt(YT) =\sum_{b\in B_n} \mu_b\epsilon_b\in P_n$, where $\mu_b$
 is the number of occurrences of $b$ in $YT$. In fact, we have
 
\begin{thm}\cite{KN}
For $\lambda=(\lambda_1,...,\lambda_n)\in \mathcal{P}_n$, $B_n(\lambda)$ is a $\frak{gl}_n$-crystal, that is isomorphism to the crystal graph of the irreducible highest weight module of $U_q(\frak{gl}_n)$ with the highest weight  $\sum_{b\in B_n} \lambda_b\epsilon_b\in P_n$.   
\end{thm}

\subsection{The Robinson-Schensted correspondence}\label{RScorr}
The Robinson–Schensted correspondence and the Littlewood–Richardson rule are well-known combinatorial algorithms that deal with the decompositions of finite-dimensional $\frak{gl}_n$
modules. 
These algorithms can be put in the framework of crystals, and yield morphisms of crystals.
Let us first recall the Schensted (row) bumping algorithm for semistandard
 tableaux.
\begin{defi}
For $\lambda\in\mathcal{P}_n$ and $YT\in B_n(\lambda)$, we define the insertion $YT\leftarrow x $ (here $x\in B_n=\{1,...,n\}$) to be the semistandard
 tableaux obtained from $YT$ by applying the following procedure:
\begin{itemize}
    \item[(I1)] if $y\le x$ for all letters $y$ on the first row, add $x$ to the right of the first row; 

    \item[(I2)] otherwise pick the leftmost node that carries a letter $x_2> x$. Replace $x_2$ on that node by $x$;

    \item[(I3)] do the same for $x_2$ and the second row, then proceed to the third row, and so on.
\end{itemize}
\end{defi}
The result of the bumping algorithm is always a semistandard tableaux.
Let us define its inverse operation as follows.
\begin{defi}\label{delop}
We denote by $\lambda\stackrel{a}{\longleftarrow} \mu$, if the node is added at the $a$-th row. That is $\mu=(\mu_1,...,\mu_n)$ with $\mu_a=\lambda_a+1$, and $\mu_i=\lambda_i$ for $i\ne a$. 
Let $YT\in B_n(\mu)$ be a semistandard tableaux. We define a new semistandard tableaux $\leftarrow_a YT\in B_n(\lambda)$, called the deletion:
\begin{itemize}
    \item[(D1)] remove the letter $x_a$ on the rightend of the $a$-th row; 

    \item[(D2)] on the $(a-1)$-th row, pick the rightmost node that carries a letter $x_{a-1}\ge x$. Replace $x_{a-1}$ on that node by $x_a$

    \item[(D3)] repeat $(D2)$ for the $(a-2)$-th row, and so on.
\end{itemize}
\end{defi}
The insertion and deletion operators are invertible, that is if $T= T'\leftarrow a$ for some $T'\in B_n(\lambda)$ and $a\in B_n$, then $T'=\leftarrow_a T$. Conversely if $a$ is the 
letter thrown away from $T$, i.e., $T'=\leftarrow_a T$, then $T=T'\leftarrow a$.

Next, associated to any word $\omega=a_1\cdots a_N\in \mathcal{W}_n$, 
we introduce a semistandard tableau $P(\omega)$ (by inserting the $N$ letters of $\omega$ into the empty tableau) and a standarde tableau $Q(\omega)$ respectively: 
\begin{itemize}
    \item we define a semistandard tableau $P(\omega)$ by 
\begin{equation}\label{Pomega}
P(\omega)=(\cdots ((\emptyset\leftarrow a_1)\leftarrow a_2)\leftarrow \cdots)\leftarrow a_N,    
\end{equation}
where $\emptyset$ is the empty tableau.

\item  We define $Q(\omega)$ to be the standard Young tableau (known as the recording tableau) with entries $1,...,N$, where $k$ is placed in the box added at the $k$-th step in the construction of $P(\omega)$. 
 
\end{itemize}

\begin{defi}\label{RSc}
The Robinson-Schensted correspondence is the bijiection 
\[RS: \ B_n^{\otimes N}\rightarrow \cup_\lambda B_n(\lambda)\times SY_N(\lambda)~:~ \omega\mapsto (P(\omega), Q(\omega)),\]
between the two sets, where the union is taken over all partitions $\lambda\in \mathcal{P}_n$ with length $N$, and $ST_N(\lambda)$ denotes the set of all standard tableaux of shape $\lambda$ with entries in $\{1,...,N\}$.
\end{defi}

Each connected component of the $\frak{gl}_n$-crystal $B_n^{\otimes N}$ is isomorphic to the crystal $B_n(\lambda)$ for some $\lambda$. One can check that for $\omega\in B_n^{\otimes N}$ and $i\in I_n$, $Q(\omega)=Q(\tilde{e}_i \omega)$, whenever $(\tilde{e}_i \omega)\ne 0$. That is the $Q$-tableaux are
 invariant under the crystal operators, and the $Q$-tableaux of $\omega$ represents the connected component of $\omega$ in $B_n^{\otimes n}$. In particular, following \cite{LT}, the map $RS$ in Definition \ref{RSc} gives an isomorphism of $\frak{gl}_n$-crystals, that gives the crystal version of the Robinson–Schensted correspondence.

\subsection{The Littlewood–Richardson rule}
Given two partitions $\mu\subset \lambda$,  a (skew) semistandard Young
tableau T of shape $\lambda/\mu$ is a filling of the diagram,
obtained by removing $\mu$ from the Young diagram of $\lambda$, with nonnegative integers
satisfying the semistandard conditions.

\begin{defi}
Given $\lambda, \mu$, and $\nu$ in $\mathcal{P}_n$ such that $\mu\subset \lambda$ and $|\lambda|=|\mu|+|\nu|$, a
semistandard tableau $YT$ of shape $\lambda/\mu $ is called a Littlewood–Richardson tableau of shape $\lambda/\mu$ with content $\nu$ if

 (1) the number of occurrences of $k$ in $YT$ is equal to $\nu_k$ for $k\ge 1$, 
 
 (2) the word $\omega(YT)$ is a lattice permutation. Here a word $\omega=\omega_1\cdots \omega_r\in \mathcal{W}_n$ is called a lattice permutation if for $1\le k\le r$ and $1\le i\le n-1$, the number of occurrences of $i$ in $\omega_1\cdots \omega_r$ is no less than that of $i+1$.
\end{defi}
 
We denote by $LR^\lambda_{\mu\nu}$ the set of all Littlewood–Richardson tableaux of shape $\lambda/\mu$ with content $\nu$.

\begin{thm}\cite{Tho}\label{LRc}
The following map is a bijection
\begin{align} 
LR:  B_n(\mu)\times B_n(\nu) &\rightarrow \cup_\lambda B_n(\lambda)\times LR^\lambda_{\mu\nu} \\
YT_1\otimes YT_2 &\mapsto \ \ (P(YT_1\otimes YT_2), Q(YT_1\otimes YT_2)),
\end{align}
where if $\omega(YT_2)=\omega_1\cdots \omega_r$, then 
\begin{equation}\label{PSY}
    P(YT_1\otimes YT_2):=((((YT_1\leftarrow\omega_1)\leftarrow \omega_2)\cdots )\leftarrow \omega_r)
\end{equation} and $Q(YT_1\otimes YT_2)$ is defined as the Littlewood–Richardson tableau of shape $\lambda/\mu$ (with $\lambda$ the shape of $P(YT_1\otimes YT_2)$) such that if $\omega_i$ is in the $k$-th row of $YT_2$ and inserted into the right hand side of \eqref{PSY} to create a node in $\lambda/\mu$, then we fill the node with $k$. 
\end{thm}

Similar to the case of the Robinson–Schensted correspondence, for
$YT_1\otimes YT_2\in B_n(\mu)\otimes B_n(\nu)$ and $i\in I_n$,
 we have $Q(YT_1\otimes YT_2)=Q(\tilde{e}_i(YT_1\otimes YT_2))$, whenever $\tilde{e}_i(YT_1\otimes YT_2)\ne 0$.
 Hence, the $Q$-tableau is invariant under the Kashiwara operators, and it represents
 the connected component of $YT_1\otimes YT_2$ in $B_n(\mu)\otimes B_n(\nu)$. In particular, 
\begin{thm}\cite{Na}\label{LRcrystal}
The map $LR$ in Theorem \ref{LRc} gives an isomorphism of $\frak{gl}_n$-crystals (that gives the crystal version of the Littlewood-Richardson rule).
\end{thm}

\subsection{The Gelfand-Tsetlin patterns}\label{GTY}
A Gelfand-Tsetlin (GT) pattern $\Lambda$ is a collection of integers $\{\lambda^{(i)}_j(\Lambda)\}_{1\le i\le j\le n-1}$ satisfying the interlacing conditions
\begin{equation}\label{interineq}
\lambda_j^{(i)}(\Lambda)-\lambda^{(i-1)}_j(\Lambda)\in\mathbb{Z}_{\ge 0}, \hspace{5mm} \lambda^{(i-1)}_j(\Lambda)-\lambda^{(i)}_{j+1}(\Lambda)\in\mathbb{Z}_{\ge 0}, \hspace{5mm} \forall \ i=1,...,n-1.
\end{equation}
Henceforth we may assume that $\lambda^{(n)}_n\ge 0$. The first row $\lambda:=\lambda^{(n)}(\Lambda)=(\lambda^{(n)}_1(\Lambda),...,\lambda^{(n)}_n(\Lambda)$
is a partition, and is called the shape of $\Lambda$. 

The set $GT_n(\lambda)$ of GT patterns of shape $\lambda$ is in one-to-one correspondence with the set $B_n(\lambda)$ of semistandard tableaux of the shape $\lambda$. Given a sequence $a = (a_1,...,a_n)$ of integers satisfying $0=a_0\le a_1\le \cdots\le a_n$, we can get another sequence $b=(b_1,...,b_k)$ by writing the letter $i$ $a_{i}-a_{i-1}$ times, for each $i$ with $a_{i-1}<a_i$, starting 
with $i = 1,2, \cdots$. Note that $b_k=a_n$. The resulting sequence is called the rearrangement of $a$. Then 
\begin{defi}\label{tau}
For any partition $\lambda$, we define $\tau$ as the bijiection
\[\tau:GT_n(\lambda)\rightarrow B_n(\lambda)~:~ \Lambda\mapsto \tau(\Lambda),\]
where 
\[\text{the $j$-th row of $\tau(\Lambda)$ is the rearrangement of sequence }(\underbrace{0,...,0}_{j},\lambda^{(j)}_j,\lambda^{(j+1)}_{j},...,\lambda^{(n)}_j).\]
\end{defi}

\section{Quantum Stokes matrices and connection matrices}\label{qbeginsection}

In this section, we first recall the quantum Stokes matrices and connection matrices of the linear equation
\begin{equation}\label{texteq}
\frac{dF}{dz}=h\left(\I u+\frac{1}{2\pi \I }\frac{\mathcal{T}}{z}\right)\cdot F,
\end{equation}
associated to a finite dimensional representation $W$ of $\frak{gl}_n$.
Following \cite{Xu, Xu2}, we then recall the regularized limits of the monodromy matrices, as well as their relations with quantum groups.

\begin{rmk}
The quantum confluent hypergeometric equation considered in \cite{Xu, Xu2, LX, TX} takes the form
\[\frac{dF}{dz}=h\left(\I u+\frac{1}{2\pi \I }\frac{T}{z}\right)\cdot F,\]
where $T=({T}_{ij})$ has entries ${T}_{ij}=e_{ij}$ for $ 1\le i,j\le n$, i.e., the residue matrix $T$ is the "transpose" of $\mathcal{T}$. All results and computations given in \cite{Xu,Xu2} can be carried over to \eqref{texteq} provided the "transpose" is taken into account. In the following, we simply list the results.
\end{rmk}

\subsection{Canonical solutions and monodromy}\label{solmono}
First, since $h$ is a positive real number, the equation \textbf{texteq} has a unique formal fundamental solution, see \cite{Xu2} for more details.
\begin{pro}\cite{Xu2}\label{formalSol}
For any nonzero real number $h$ and $u\in\h_{\rm reg}(\mathbb{R})$, the ordinary differential equation \eqref{texteq} has a unique formal fundamental solution taking the form \begin{eqnarray}\label{formalsum}
\widehat{F}(z)=\widehat{H}(z) e^{{h\I uz}}z^{h[\mathcal{T}]}, \ \ \ {\it for} \ \widehat{H}=1+H_1z^{-1}+H_2z^{-2}+\cdot\cdot\cdot, \end{eqnarray}
where each coefficient $H_m\in{\rm End}(W)\otimes{\rm End}(\mathbb{C}^n)$, and $[\mathcal{T}]$ denotes the diagonal part of $ 
 T$, i.e., $[\mathcal{T}]=\sum_{k} e_{kk}\otimes E_{kk}.$
\end{pro}
\begin{rmk}
In \cite{TX}, the universal expression of each $H_m=H_{m,ij}\otimes E_{ij}$ was given. That is explicit elements $h_{m,ij}\in U(\frak{gl}_n)$ are given, such that for any representation $\rho:U(\frak{gl}_n)\rightarrow W$, $H_{m,ij}=\rho(h_{m,ij})\in {\rm End}(W)$. Furthermore, a relation between $H_m$ and the Yangian of $\frak{gl}_{n-1}$ was given.
\end{rmk}

%The radius of convergence of the formal power series $\hat{H}(z)$ in Proposition \ref{formalSol} is in general zero. However, it follows from the general principle of differential equations with irregular singularities that (see e.g., \cite{Balser,LR}) the Borel-Laplace transform of $\hat{H}(z)$ along any admissible direction gives a holomorphic function in each Stokes supersector around $z=\infty$. In this way, one gets actual solutions of \eqref{texteq} on each Stokes supersector. These sectors are determined by the irregular term $hu$ of the differential equation as follows.

\begin{rmk}
The anti-Stokes rays of the equation \eqref{texteq} are the directions along which $e^{h\I (u_k-u_j)z}$ decays most rapidly as $z\rightarrow \infty$ for some $u_k\ne u_j$. In our case, since $h\in \mathbb{R}$ and $u\in\h_{\rm reg}(\mathbb{R})$, the anti-Stokes rays are just the two imaginary axis. Then the standard resummation procedure that (see e.g., \cite{Balser,LR}) the Borel-Laplace transform of $\hat{H}(z)$, along any direction $d$ but the anti-Stokes rays, gives a holomorphic function with the asymptotics $\hat{H}(z)$ as $z\rightarrow\infty$ within a proper sector centered on $d$.
\end{rmk}

Let us choose the branch of ${\rm log}(z)$, which is real on the positive real axis, with a cut along $i\mathbb{R}_{\ge 0}$. 
For any two real numbers $a, b$, we denote by $S(a,b):=\{z\in\mathbb{C}~|~a<{\rm arg}(z)<b\}$ an open sector with opening
angle $b-a>0$.
It follows from the general principle of differential equations that

\begin{thm}(see e.g., \cite{Xu})\label{quantumsol}
For any $u\in\h_{\rm reg}(\mathbb{R})$ and integer $k$, there
is a unique (therefore canonical) multivalued holomorphic fundamental solution $F_{k}(z;u)\in {\rm End}(W)\otimes {\rm End}(\mathbb{C}^n)$ of \eqref{texteq} such that 
\begin{equation*}
F_{k}(z;u)\cdot e^{-\I h uz}\cdot z^{\frac{-h[\mathcal{T}]}{2\pi\I }}\sim \hat{H}(z), \ \ \ \text{as $z\rightarrow\infty$ within the sector } {S}\left((k-1)\pi,(k+1)\pi\right).
\end{equation*}
\end{thm}
It follows from Watson's uniqueness theorem that $F_{k+2}(z;u)=F_{k}(ze^{2\pi\I};u)\cdot e^{h[\mathcal{T}]}$. Thus we only need to study the solutions $F_+(z;u):=F_0(z;u)$ and $F_-(z;u):=F_{-1}(z;u)$ with $k=0,-1$ respectively. These two solutions are specified by the asymptotics
\begin{align}\label{asym1}
\lim_{z\rightarrow\infty}F_{+}(z;u)\cdot e^{-\I h uz}\cdot z^{\frac{-h[\mathcal{T}]}{2\pi\I }}&=1, \ \ \ as \ \ \ z\in {S}(-\pi,\pi),
\\ \label{asym2}
\lim_{z\rightarrow\infty}F_{-}(z;u)\cdot e^{-\I h uz}\cdot z^{\frac{-h[\mathcal{T}]}{2\pi\I }}&=1, \ \ \ as \ \ \ z\in {S}(-2\pi,0).
\end{align}
In the rest of the paper, let us always assume 
\begin{equation}\label{uUid}
    u={\rm diag}(u_1,...,u_n)\in U_{\rm id}\subset \h_{\rm reg}(\mathbb{R}), \text{ i.e., } u_1<u_2<\cdots <u_n.
\end{equation}
\begin{defi}\cite{Xu, Xu2}\label{qStokes}
The {\it quantum Stokes matrices} of the differential equation \eqref{texteq}, associated to the given representation $W$, are the elements $S_{\pm}(u)^W\in {\rm End}(W)\otimes {\rm End}(\mathbb{C}^n)$ determined by
\begin{align*}
F_{+}(z)&=F_{-}(z)\cdot S_{+}(u)^W, \\
F_{-}(z)&=F_{+}(ze^{-2\pi\I })\cdot e^{-h[\mathcal{T}]}S_{-}(u)^W.
\end{align*}
\end{defi}
\begin{lem}\label{upper}
The matrices $S_{\pm}(u)$ are triangular with $1\in {\rm End}(L(\lambda))$ along the diagonal.
\end{lem}
\begin{proof} By Definition \ref{qStokes}, we have
\[\left(F_-(z)e^{-\I h uz} z^{\frac{-h[\mathcal{T}]}{2\pi\I }}\right)^{-1} \left(F_+(z) e^{-\I h uz} z^{\frac{-h[\mathcal{T}]}{2\pi\I }}\right)  =z^{\frac{h[\mathcal{T}]}{2\pi\I }} e^{\I h uz} S_{+}(u)^W e^{-\I h uz}z^{\frac{-h[\mathcal{T}]}{2\pi\I }}.\]
Then using the asymptotics \eqref{asym1} and \eqref{asym2}, we have
\[z^{\frac{h[\mathcal{T}]}{2\pi\I }} e^{\I h uz} S_{+}(u)^W e^{-\I h uz}z^{\frac{-h[\mathcal{T}]}{2\pi\I }}\rightarrow 1, \] 
as $z\rightarrow \infty \text{ within } \Sect_+\cap \Sect_-.$
Let us write $e^{\frac{h[\mathcal{T}]}{2}}S_{+}(u)=\sum_{ij} S_{+}(u)_{ij}\otimes E_{ij}\in {\rm End}(W)\otimes {\rm End}(\mathbb{C}^n)$. Since the exponential $e^{\I huz}$ dominate, we must have, for any $i\ne j$, $e^{\I h(u_i-u_j)z}S_{+}(u)_{ij}\otimes E_{ij} \rightarrow 0.$
It implies that  $( S_{+}(u))_{ij}=\delta_{ij}$ unless $e^{\I z h(u_i-u_j)}\rightarrow 0$ as $z\rightarrow \infty$ within ${S}(-\pi,\pi)\cap {S}(-2\pi,0)$. Thus due to the assumption \eqref{uUid}, $S_{+}(u)$ is upper triangular with diagonal entries equal to $1$. The argument for $S_{-}(u)$ is the same once the change in branches of ${\rm log}(z)$ is taken into account.
\end{proof} 

Let us now introduce a canonical solution around $z=0$.
Set $m={\rm dim}(W)$, then $\mathcal{T}\in {\rm End}(W)\otimes {\rm End}(\mathbb{C}^n)$ can be seen as $mn\times mn$ matrix. Following Lemma \ref{diaglem} the difference between any two eigenvalues of the $mn\times mn$ matrix $\mathcal{T}$ are integers. Thus the $mn\times mn$ linear system \eqref{texteq} is non-resonant, i.e., the difference between any two eigenvalues of the residue matrix $\frac{h T}{2\pi\I}$ of the $mn\times mn$ linear system \eqref{texteq} does not take non-zero integers (actually they are purely imaginary numbers of the form $\frac{h}{2\pi\I}\times c$ with some integer $c$ and $h\in\mathbb{R}$). Therefore, following a general theory of linear system of equations (see e.g \cite[Chapter 2]{Balser}), we have
\begin{lem}\label{le:nr dkz}
There is a multivalued holomorphic solution
$F_0(z)$ of the differential equation \eqref{introqeq}, such that 
\[F_0\cdot z^{-\frac{h \mathcal{T}}{2\pi\I }}\rightarrow 1 \ \text{ as } z\rightarrow 0.\]
\end{lem}

\begin{defi}\label{connectionmatrix}
The {\it quantum connection 
matrix} $C(u)^W\in {\rm End}(W)\otimes {\rm End}(\mathbb{C}^n)$ of the system \eqref{texteq}, associated to the given representation $W$, is determined by 
\[F_0(z)=F_+(z)\cdot C(u)^W. \]
\end{defi}

The connection matrix is related to the Stokes matrices
by the following monodromy relation, which follows from the fact that a simple negative loop (i.e., in clockwise direction) around $0$ is a simple positive loop (i.e., in anticlockwise direction) around $\infty$: 
\begin{eqnarray}\label{monodromyrelation}
C(u)e^{-h\mathcal{T}}C(u)^{-1}=e^{-h[\mathcal{T}]}S_-(u) S_+(u)\in {\rm End}(W)\otimes {\rm End}(\mathbb{C}^n).
\end{eqnarray}

\subsection{Quantum Stokes matrices and connection matrices at a caterpillar point}\label{sec:qStokescat} 

Recall that we denote by $u_{\rm cat}$ the limit as in \eqref{limit}. The regularized limits of the quantum Stokes matrices $S_{\pm}(u)$ as $u\rightarrow u_{\rm cat}$ from $u\in U_{\rm id}$ are defined as follows.

\begin{defi}\label{qScat}
The quantum Stokes matrices at $u_{\rm cat}$ (with respect to the choice of $U_{\rm id}$) are the upper and lower $n\times n$ triangular matrices $S_{\pm}(u_{\rm cat})^W$ having the same diagonal part, with entries valued in ${\rm End}(W)$, such that (the blocked Gauss decomposition)
\[S_{-}(u_{\rm cat})^WS_{+}(u_{\rm cat})^W\]
equals to the limit of 
\begin{equation}\label{gauge}
{\rm Ad}{\left((\frac{1}{u_{2}-u_{1}})^{\frac{{\rm log}(\delta_1(S^W_{-})\delta_1(S^W_{+}))}{2\pi\I }}\overrightarrow{\underset{k=2,...,m-1}{\prod} }(\frac{u_{k}-u_{k-1}}{u_{k+1}-u_{k}})^{\frac{{\rm log}(\delta_k(S^W_{-})\delta_k(S^W_{+}))}{2\pi\I }}\right)}  \left(S_{-}(u)^W S_{+}(u)^W\right),
\end{equation}
as $u\rightarrow u_{\rm cat}$ from $U_{\rm id}$, i.e., $\frac{u_{k+1}-u_{k}}{u_{k}-u_{k-1}}\rightarrow +\infty$ for all $k=2,...,m-1$ and $u_2-u_1\rightarrow 0$,
where $\delta_k(S_\pm(u))\in {\rm End}(W)\otimes {\rm End}(\mathbb{C}^n)$ is the matrix with entries
\[\delta_k(S_\pm(u))_{ij}:=\left\{
          \begin{array}{lr}
            (S_\pm(u))_{ij},   & \text{if} \ \ 1\le i, \ j \le k  \ \text{ or } \ i=j\\
           0, & \text{otherwise},
             \end{array}
\right. \]

\end{defi}

\begin{defi}\label{qCcat}
The quantum connection matrix $C(u_{\rm cat})^W\in {\rm End}(W)\otimes {\rm End}(\mathbb{C}^n)$ at $u_{\rm cat}$ (with respect to the choice of $U_{\rm id}$) is the limit
\begin{equation}\label{Cgauge}
C(u_{\rm cat})^W=\lim_{u\rightarrow u_{\rm cat}} {\rm Ad}{\left(\left(\frac{1}{u_{2}-u_{1}}\right)^{\frac{{\rm log}(\delta_1(S^W_{-})\delta_1(S^W_{+}))}{2\pi\I }}\overrightarrow{\underset{k=2,...,m-1}{\prod} }\left(\frac{u_{k}-u_{k-1}}{u_{k+1}-u_{k}}\right)^{\frac{{\rm log}(\delta_k(S^W_{-})\delta_k(S^W_{+}))}{2\pi\I }}\right)}  C(u)^W.
\end{equation}
\end{defi}
See Remark \ref{rmkspin} for the meaning of the regularization term (the term under the Adjoint action).
We remark that the monodromy relation is preserved, that is
\begin{eqnarray}\label{qcatmono}
S_{-}(u_{\rm cat})S_{+}(u_{\rm cat})=C(u_{\rm cat})e^{-h\mathcal{T}} C(u_{\rm cat})^{-1}.
\end{eqnarray}

In \cite{Xu}, we have used the method of isomonodromy deformation to prove that, as in the limit $\frac{u_j-u_{j-1}}{u_{j-1}-u_{j-2}}\rightarrow +\infty$ for all $j=3,...,n$, the differential equation \eqref{texteq} is "decoupled into" a set of equations for functions $F_k(z)\in {\rm End}(W)\otimes{\rm End}(\mathbb{C}^n)$ with $k=2,...,n$,
\begin{equation}\label{eqk}
\frac{dF_k}{dz}=h\left(\I E_{kk}+\frac{1}{2\pi \I }\frac{\delta_k(\mathcal{T})}{z}\right)\cdot F_k,
\end{equation}
where $E_{kk}={\rm diag}(0,,,0,1,0,...,0)$ and \[\delta_k(\mathcal{T})_{ij}=\left\{
          \begin{array}{lr}
             \mathcal{T}_{ij},   & \text{if} \ \ 1\le i, j\le k, \ \text{or} \ i=j  \\
           0, & \text{otherwise}.
             \end{array}
\right.\]
More precisely, in the limit, the Stokes and connection matrices of \eqref{texteq} are given by explicit combinations of the monodromy of equations \eqref{eqk} for all $k=2,...,n$ (up to the regularization terms as in \eqref{gauge} and \eqref{Cgauge}).
\begin{thm}\cite{Xu}
The connection matrix at the caterpillar point $u_{\rm cat}$ is given by 
\begin{equation}\label{leftprod}
    C(u_{\rm cat})^{W}=\overrightarrow{\prod_{k=2,...,n}}C^{(k)}\in {\rm End}(W)\otimes {\rm End}(\mathbb{C}^n),
\end{equation}
where $C^{(k)}\in {\rm End}(W)\otimes {\rm End}(\mathbb{C}^n)$ is the connection matrix of equation \eqref{eqk}, and the product $\overrightarrow{\prod}$ is taken with the index $i$ to the right of $j$ if $i>j$. The Stokes matrices at the limit $u_{\rm cat}$ is then given by the Gauss decomposition \eqref{qcatmono}.
\end{thm}

The differential equation \eqref{eqk} is exactly solvable by the confluent hypergeometric functions, using which one can express its connection matrix $C^{(k)}$ explicitly via gamma functions (as in Theorem \ref{qexStokes}), see \cite{LX, Xu} for more details. In this way, we obtain the explicit expression of the Stokes and connection matrices at the limit point $u_{\rm cat}$. In \cite{Xu}, the explicit expression are given under the Gelfand-Tsetlin basis of the associated representation.

\subsection{Quantum minors and Gelfand-Tsetlin basis}\label{GZbasis}
To write down the explicit expression of $S_{\pm}(u_{\rm cat})$ associated to an irreducible representation $L(\lambda)$, let us first introduce the quantum minors of $\mathcal{T}$ and their actions on the Gelfand-Testlin basis. 

\subsubsection*{Quantum minors and diagonalization of $\mathcal{T}$}
The matrix $\mathcal{T}(x)=\mathcal{T}-x \operatorname{Id}$ is the characteristic matrix of the matrix $\mathcal{T}=(e_{ji})_{i,j=1,...,n}$, where $x$ is an indeterminate commuting with all generators $e_{ji}$, for $1\leq i,j \leq n$.

\begin{defi}
For $1\leq m \leq n$, given two sequences, $a=\{a_1,\cdots,a_m\}$ and $b=\{b_1,\cdots,b_m\}$, with elements in $\{1,2,\cdots,n-1,n\}$, the corresponding quantum
minor of the matrix $\mathcal{T}(x)$ is defined by 
\begin{align*}
\Delta^{a_1,...,a_m}_{b_1,...,b_m}(\mathcal{T}(x)):=\sum_{\sigma\in S_m}(-1)^\sigma \mathcal{T}(x_1)_{a_{\sigma(1)},b_1}\cdots \mathcal{T}(x_m)_{a_{\sigma(m)},b_m}\in U({\frak {gl}}_n)[x],
\end{align*}
where $x_k:=x+k-1$ for $k=1,...,m$, and $(-1)^{\sigma}$ means the signature of the permutation $\sigma$ in the symmetry group $S_m$ of $m$ elements.
\end{defi}
Suppose that $\Delta^{a_1,...,a_m}_{b_1,...,b_m}\left(\mathcal{T}\left(x\right)\right)=\sum_{i=0}^{n-1}r_i x^i$ with coefficient $r_i\in U\left({\mathfrak {gl}}_n\right)$. For any element $\zeta$ (that may not commute with $e_{ij}$ for all $1\le i,j\le n$), we define $\Delta^{a_1,...,a_m}_{b_1,...,b_m}\big(\mathcal{T}(\zeta)\big)$ in the following convention
\begin{equation*}
\Delta^{a_1,...,a_m}_{b_1,...,b_m}\big(\mathcal{T}(\zeta)\big)=\zeta^{n-1} r_{n-1}+\zeta^{n-2} r_{n-2}+\cdots +\zeta r_{1}+r_0.
\end{equation*}

\begin{defi}\label{qroots}
For any $1\le k\le n$, let 
$
M_k(\zeta):=\Delta^{1,...,k}_{1,...,k}\left(T(\zeta)\right)$
be the upper left $k\times k$ quantum-minor of $T$. 
Let $\{\zeta^{(k)}_i\}_{1\le i\le k\le n}$ denote the roots of $M_k(\zeta) = 0$ for all $k=1,...,n$ (in an appropriate splitting extension).
\end{defi}

It implies that 
\begin{lem}\label{diaglem}
The matrix $\mathcal{T}$ can be diagonalized to 
\[D=\begin{pmatrix}
    \zeta^{(n)}_1 & & & \\
    & \ddots & \\
    & &  \zeta^{(n)}_n 
    
    \end{pmatrix}.\]
\end{lem}
In fact, one can explicitly write down the conjugation that diagonalizes $\mathcal{T}$. Let $P$ and $Q$ be $n\times n$ matrices defined by
\begin{align}\label{P}
&P=\left(\begin{array}{ccc}
\left(-1\right)^{1+1} \Delta^{\hat{1} ,..., n-1 }_{1 ,..., n-1}\left(\mathcal{T}\left(\zeta^{\left(n\right)}_{1}+1\right)\right) & \cdots & \left(-1\right)^{1+n} \Delta^{\hat{1} ,..., n}_{1 ,..., n}\left(\mathcal{T}\left(\zeta^{\left(n\right)}_{n}+1\right)\right) \\ 
\vdots &  \ddots & \vdots \\
\left(-1\right)^{n+1} \Delta^{1 ,..., \widehat{n} }_{1 ,..., n-1}\left(\mathcal{T}\left(\zeta^{\left(n\right)}_{1}+1\right)\right) & \cdots & \left(-1\right)^{n+n} \Delta^{1 ,..., \widehat{n} }_{1 ,..., n-1}\left(\mathcal{T}\left(\zeta^{\left(n\right)}_{n}+1\right)\right)
\end{array}\right)\cdot \frac{1}{\mathcal{N}},\\ \label{Q}
&Q=\frac{1}{\mathcal{N}}\cdot \left(\begin{array}{ccc}
\left(-1\right)^{1+1} \Delta^{1 ,..., n-1}_{\hat{1} ,..., n}\left(\mathcal{T}\left(\zeta^{\left(n\right)}_1+1\right)\right) & \cdots & \left(-1\right)^{1+n} \Delta_{1 ,..., \widehat{n} }^{1 ,..., n-1}\left(\mathcal{T}\left(\zeta^{\left(n\right)}_1+1\right)\right) \\
\vdots & \ddots & \vdots \\
\left(-1\right)^{n+1} \Delta^{1 ,..., n-1}_{\hat{1} ,..., n}\left(\mathcal{T}\left(\zeta^{\left(n\right)}_{n}+1\right)\right) & \cdots & \left(-1\right)^{n+n} \Delta_{1 ,..., \widehat{n}}^{1 ,..., n-1}\left(\mathcal{T}\left(\zeta^{\left(n\right)}_{n}+1\right)\right)
\end{array}\right).
\end{align}
where $\mathcal{N}={\rm diag}\left(\mathcal{N}_1,....,\mathcal{N}_{n}\right)$ is a diagonal matrix with the diagonal elements 
\[\mathcal{N}_i=\sqrt{\prod_{l=1,l\ne i}^{n-1}\left(\zeta^{\left(n\right)}_l-\zeta^{\left(n\right)}_i\right)\prod_{l=1}^{n-1}\left(\zeta^{\left(n-1\right)}_l-\zeta^{\left(n\right)}_i-1\right)}.\]
Then one can check that 
\begin{lem}\label{diagT}[See e.g., \cite{LX}]
The matrices satisfy $PQ=1=QP$ and $Q\mathcal{T}P=D$.
\end{lem}

\subsubsection*{Gelfand-Tsetlin basis}
It is known that the subalgebra, generated in $U({\frak {gl}}_n)$ by the coefficients in all minors $M_k(\zeta)$ (given in Definition \ref{qroots}) for $1\le k\le n$, is a maximal commutative subalgebra and is called the Gelfand-Tsetlin subalgebra. 
Let us now take a highest weight representation $L(\lambda)$. The action of the Gelfand-Tsetlin subalgebra on $L(\lambda)$ has simple spectrum and one corresponding set $E_{GT}(\lambda)$ of orthonormal basis $\xi_\Lambda$ of $L(\lambda)$ is called the Gelfand-Testlin basis, with basis vectors $\xi_\Lambda$ parameterized by all patterns $\Lambda\in GT_n(\lambda)$.

Then the roots $\{\zeta^{(k)}_i\}_{1\le i\le k}$ of $M_k(\zeta) = 0$ can be specified by explicit actions on the basis. The action of quantum minors can also be explicitly expressed under the basis. The following proposition can be found in \cite[Chapter 2]{Molev} (provided that the normalization of the orthogonal basis in \cite{Molev} is taken into account).

\begin{pro}\label{GZaction}
For $1\le i\le k\le n$, the actions of $e_{kk}, \zeta^{(k)}_i, \Delta^{1,...,k}_{1,...,k-1,k+1}\left(\mathcal{T}(\zeta^{(k)}_i)+1\right)\in {\rm End}(L(\lambda))$ on the basis ${\xi_\Lambda}$ of $L(\lambda)$ are given by
\begin{align}
e_{kk}\cdot \xi_\Lambda &=\left(\sum_{i=1}^k\lambda^{(k)}_i(\Lambda)-\sum_{i=1}^{k-1}\lambda^{(k-1)}_i(\Lambda)\right)\xi_\Lambda,\\ \label{eigenvalues}
\zeta^{(k)}_i\cdot \xi_\Lambda &=\left(\lambda^{(k)}_i(\Lambda)-i+1\right)\xi_\Lambda ,
\end{align}
and 
\begin{align}\nonumber
&\Delta^{1,...,k}_{1,...,k-1,k+1}\left(\mathcal{T}(\zeta^{(k)}_i+1)\right) \cdot \xi_\Lambda \\ \label{ashift}
=&(-1)^{k+i}\sqrt{-\frac{\prod_{l=1,l\ne i}^{k}(\zeta^{(k)}_i-\zeta^{(k)}_l)\prod_{l=1}^{k+1}(\zeta^{(k)}_i-\zeta^{(k+1)}_l)\prod_{l=1}^{k-1}(\zeta^{(k)}_i-\zeta^{(k-1)}_l+1)}{\prod_{l=1,l\ne i}^{k}(\zeta^{(k)}_i-\zeta^{(k)}_l+1)}}\cdot \xi_{\Lambda-\delta^{k}_i} 
\end{align}
where the pattern $\Lambda-\delta^{(k)}_i$ is obtained from $\Lambda$ by replacing $\lambda^{(k)}_i$ by $\lambda^{(k)}_i-1$. It is supposed that $\xi_{\Lambda-\delta^{(k)}_i}$ is zero if $\Lambda-\delta^{(k)}_i$ is not a pattern.
\end{pro}

\subsection{The explicit expression of quantum Stokes matrices and connection matrix at a caterpillar point}\label{proofofqS}

\begin{thm}\cite{Xu}\label{qexStokes}
The connection matrix 
\begin{equation}
    C(u_{\rm cat})^{L(\lambda)}=\overrightarrow{\prod_{k=2,...,n}}C^{(k)}=\left(\overrightarrow{\prod_{k=2,...,n}}\widetilde{C^{(k)}}\right)\cdot P\in {\rm End}(L(\lambda))\otimes {\rm End}(\mathbb{C}^n),
\end{equation}
where the entries of the matrix $\widetilde{C^{(k+1)}}$ are 

(1). for $1\le j\le k+1, 1\le i\le k$,
\begin{align*} \widetilde{C_{ij}^{(k+1)}}=&\frac{-(\I h)^{\frac{h(\zeta^{(k)}_i-\zeta^{(k+1)}_j)}{2\pi\I }}e^{\frac{h(\zeta^{(k+1)}_j-2\zeta^{(k)}_i)}{2}}\cdot\sqrt{\prod_{l=1}^{k}(\zeta^{(k+1)}_j-\zeta^{(k)}_l-1)}}{\sqrt{\prod_{l=1, l\ne j}^{k+1}(\zeta^{(k+1)}_j-\zeta^{(k+1)}_l) \prod_{l=1,l\ne i}^{k-1}\left(\zeta^{\left(k\right)}_l-\zeta^{\left(k\right)}_i\right)\prod_{l=1}^{k-1}\left(\zeta^{\left(k-1\right)}_l-\zeta^{\left(k\right)}_i-1\right)}} \\
&\times \frac{\prod\limits_{l=1}^{k+1}\Gamma(1+h\frac{\zeta^{(k+1)}_j-\zeta^{(k+1)}_l}{2\pi \I })\prod\limits_{l=1}^{k}\Gamma(1+h\frac{\zeta^{(k)}_i-\zeta^{(k)}_l-1}{2\pi \I })}{\prod\limits_{l=1,l\ne i}^{k}\Gamma(1+h\frac{\zeta^{(k+1)}_j-\zeta^{(k)}_l-1}{2\pi \I })\prod\limits_{l=1,l\ne j}^{k+1}\Gamma(1+h\frac{\zeta^{(k+1)}_l-\zeta^{(k)}_i}{2\pi \I })}\cdot \frac{1}{2\pi\I} \Delta^{1,...,k}_{1,...,k-1,k+1}\left(\mathcal{T}(\zeta^{(k)}_i)\right);
\end{align*}
(2). for $\ 1\le j\le k+1$,
\begin{eqnarray*} \widetilde{C_{k+1,j}^{(k+1)}}=(\I h)^{\frac{h(e_{k+1,k+1}-\zeta^{(k+1)}_j-k)}{2\pi\I }}e^{\frac{-he_{k+1,k+1}}{2}} \cdot \sqrt{\frac{\prod_{v=1}^{k}(\zeta^{(k+1)}_j-\zeta^{(k)}_v+1)}{\prod_{v=1, v\ne j}^{k+1}(\zeta^{(k+1)}_j-\zeta^{(k+1)}_v) }}\frac{\prod_{l=1}^{k+1}\Gamma(1+h\frac{\zeta^{(k+1)}_j-\zeta^{(k+1)}_l}{2\pi \I })}{\prod_{l=1}^{k}\Gamma(1+h\frac{\zeta^{(k+1)}_j-\zeta^{(k)}_l+1}{2\pi \I })};
\end{eqnarray*}
(3). \[\ \widetilde{C_{ii}^{(k+1)}}=1  \ \ \text{for} \ \ k+2\le i\le n; \hspace{5mm} \widetilde{C_{ij}^{(k+1)}}=0 \ \ \text{otherwise}.\]
\end{thm}

\begin{rmk}
    
In \cite{Xu}, the matrix $\widetilde{C^{(k)}}$ was called the normalized connection matrix of the equation \eqref{eqk}. It satisfies $\widetilde{C^{(k)}}=(P_{k-1})^{-1}{C^{(k)}} P_k$. Here $P_k=\left(\begin{array}{cc}
    P^{(k)} & 0  \\
    0 & {\rm Id}_{n-k}
  \end{array}\right)$, and $ P^{(k)}$ is the $k\times k$ matrix, that diagonalizes the upper left $k\times k$ submatrix of $\mathcal{T}$, as defined in \eqref{P} for $n=k$. The connection matrix of \eqref{eqk} was calculated in \cite{LX} (provided that the "transpose" of $\mathcal{T}$ is taken into account).
\end{rmk}
Using the monodromy relation (Gauss decomposition) \eqref{qcatmono}, we have 
\begin{thm}\cite{Xu}\label{explicitS}
For any $1\le k\le n-1$, the $(k,k+1)$-entry of $S_{+}(u_{\rm cat})^{L(\lambda)}$, as an element in ${\rm End}(L(\lambda))$ is given by
\begin{align*}
&(S_{+}(u_{\rm cat})^{L(\lambda)}_{k,k+1}=
 \left(\frac{h}{2\pi\I }\right)^k (\I h)^{\frac{h\small{(e_{kk}-e_{k+1,k+1}-1})}{2\pi \I }} \\
&\times  \sum_{i=1}^k\left(\frac{\prod_{l=1,l\ne i}^{k}\Gamma\left(h\frac{\zeta^{(k)}_i-\zeta^{(k)}_l}{2\pi \I }\right)}{\prod_{l=1}^{k+1}\Gamma\left(1+h\frac{\zeta^{(k)}_i-\zeta^{(k+1)}_l-1}{2\pi \I }\right)}\frac{\prod_{l=1,l\ne i}^{k}\Gamma\left(1+h\frac{\zeta^{(k)}_i-\zeta^{(k)}_l-1}{2\pi \I }\right)}{\prod_{l=1}^{k-1}\Gamma\left(1+h\frac{\zeta^{(k)}_i-\zeta^{(k-1)}_l}{2\pi \I }\right)}\cdot \Delta^{1,...,k}_{1,...,k-1,k+1}\left(\mathcal{T}(\zeta^{(k)}_i+1)\right)\right).\end{align*}
The other entries of $S_{+}(u_{\rm cat})^{L(\lambda)}$ are also given by explicit formula, and $e^{-h[\mathcal{T}]}S_{-}(u_{\rm cat})^{L(\lambda)}e^{h[\mathcal{T}]}$ is the complex conjugation of $S_{+}(u_{\rm cat})^{L(\lambda)}$. 
\end{thm}

\begin{rmk}\label{rmkspin}
The regularization term (the term under the Adjoint action) in Definition \ref{qScat} can now be better understood. It was proved in \cite{Xu} that 
\begin{equation}\label{hnonuni}
S_{+}(u)_{k,k+1}={l}^{(+)}_{k,k+1}(u)+o(l^{(\pm)}_{k,k+1}(u)), \hspace{3mm} \text{as} \ u\rightarrow u_{\rm cat} \ \text{ from } U_{\rm id},
\end{equation}
where the leading term ${l}^{(+)}_{k,k+1}(u)$ is given by the expression $S_{+}(u_{\rm cat})_{k,k+1}$ in Theorem \ref{explicitS} provided we replace $\Delta^{1,...,k}_{1,...,k-1,k+1}\left(\mathcal{T}(\zeta^{(k)}_i+1)\right)$ by 
\begin{equation}\label{modify0}
(u_k-u_{k-1})^{\frac{he_{kk}-h\zeta^{(k)}_i}{2\pi\I }}\left({u_{k+1}-u_k}\right)^{\frac{h\zeta^{(k)}_i-he_{k+1,k+1}}{2\pi\I }}\cdot \Delta^{1,...,k}_{1,...,k-1,k+1}\left(\mathcal{T}(\zeta^{(k)}_i+1)\right).
\end{equation}
Following Proposition \ref{GZaction}, under the Gelfand-Tsetlin basis, the entry $S_{+}(u)_{k,k+1}\in {\rm End}(L(\lambda))$ satisfies 
\begin{equation}\label{fastspin}
S_{+}(u)_{k,k+1}\cdot \xi_\Lambda= \sum_{j=1}^k g_j(\Lambda;h)\times e^{\I h \theta_j(u)}\xi_{\Lambda-\delta^{(k)}_j}+\sum_{k=2}^{n-1}\mathcal{O}\left(\frac{u_{k}-u_{k-1}}{u_{k+1}-u_{k}}\right),    
\end{equation}
as $u\rightarrow u_{\rm cat}$ from $U_{\rm id}$,
where $g_j(\Lambda;h)$ are complex functions and $\theta_j(u)$ are certain real functions (here each $\mathcal{O}\left(\frac{u_{k}-u_{k-1}}{u_{k+1}-u_{k}}\right)$ stands for a remainder whose vector norm is less than $M \times\left(\frac{u_{k+1}-u_{k}}{u_{k}-u_{k-1}}\right)^{-1}$ for a positive real number $M$ as $\frac{u_{k+1}-u_{k}}{u_{k}-u_{k-1}}$ big enough). 

From \eqref{fastspin}, we see that it is the fast spin terms $e^{\I h \theta_j(u)}$ that prevent the Stokes matrix entries from having a limit. Thus, the regularization term in Definition \ref{qScat} plays the role of gauging the fast spin terms $e^{\I h \theta_j(u)}$ out. 
  
\end{rmk}
\subsection{Quantum groups and coproduct of quantum Stokes matrices}\label{qgStokes}

\begin{defi}\label{defqg}
The quantum group $U_q(\frak{gl}_n)$ is a unital associative algebra with generators $q^{\pm h_i}, e_j, f_j,$ $1\le j\le n-1, 1\le i\le n$ and the relations:
\begin{itemize}
    \item for each $1\le i\le n$, $1\le j\le n-1$,
\[q^{h_i}q^{-h_i}=q^{-h_i}q^{h_i}=1, \
q^{h_i}e_jq^{-h_i}=q^{\delta_{ij}}q^{-\delta_{i,j+1}}e_j, \ q^{h_i}f_jq^{-h_i}=q^{-\delta_{ij}}q^{\delta_{i,j+1}}f_j;
\]
\item for each $1\le i,j\le n-1$,
\[
[e_i,f_j] = \delta_{ij} \frac{q^{h_i-h_{i+1}}-q^{-h_i+h_{i+1}}}{q-q^{-1}};
\]
\item for $|i-j|=1$, 
\[
e_i^2e_j - (q+q^{-1})e_ie_je_i + e_je_i^2=0, 
\]
\[
f_i^2f_j - (q+q^{-1})f_if_jf_i + f_jf_i^2=0,
\]
and for $|i-j|\ne 1$, $[e_i,e_j]=0=[f_i,f_j]$.
\end{itemize}
Furthermore, $U_q(\frak{gl}_n)$ is a Hopf algebra with the comultiplication $\Delta$
\begin{align}
\Delta(q^{h_i})=q^{h_i}\otimes q^{h_i}, \ \ \Delta(e_i)=1\otimes {e_i}+e_i\otimes q^{h_i-h_{i+1}}, \ \ \Delta({f_i})=q^{-h_i+h_{i+1}}\otimes {f_i}+f_i\otimes 1. 
\end{align}
\end{defi}

Theorem \ref{quantumgpucat} gives the realization of (the associative algebra structure of) quantum groups via the Stokes matrices, it is natural to mimic the coproduct in quantum groups by defining the operators: 
\begin{defi}
For any two finite dimensional representations $W$ and $K$, we define
\begin{align}\label{coS}
    \Delta(S_{+}(u)_{i,i+1})^{W\otimes K}:&=q^{e_{ii}-e_{i+1,i+1}} \otimes S_{+}(u)^K_{i,i+1}+S_{+}(u)^W_{i,i+1}\otimes 1 \in {\rm End}(W\otimes K),\\
    \Delta(S_{-}(u)_{i+1,i})^{W\otimes K}:&=q^{e_{i+1,i+1}-e_{ii}} \otimes S_{-}(u)^K_{i+1,i}+S_{-}(u)^W_{i+1,i}\otimes q^{2e_{i+1,i+1}-2e_{ii}} \in {\rm End}(W\otimes K).
\end{align}
\end{defi}

In particular, using Theorem \ref{quantumgpucat} associated to representations $L(\mu)$ and $L(\nu)$, there are two representations of the quantum group $U_q(\frak{gl}_n)$ on the tensor space: the first is via the Stokes matrices $S_{\pm}(u)^{L(\mu)\otimes L(\nu)}$ associated to the tensor representation, and the second is via the coproduct of the Stokes matrices $\Delta(S_{\pm}(u))^{L(\mu)\otimes L(\nu)}$.
Note that the difference of the two representations measures the difference of the coproduct on the quantum group $U(\frak{gl}_n)$ and the coproduct on the undeformed universal enveloping algebra $U(\frak{gl}_n).$

In general, the Stokes and connection matrices are (new) highly transcendental matrix functions depending on the chosen representation $W$ and the variables $u={\rm diag}(u_1,...,u_n)$. For fixed $W$ they are real analytic functions of $u\in \h_{\rm reg}(\mathbb{R})$. However, if we want to study the dependence of the Stokes matrices on the representation spaces $W$, we can take a special $u$ to minimize analytical complexity as much as possible. Our strategy is to take the reference point $u$ as the limit $u_{\rm cat}$. The quantum group relation, satisfied by the Stokes matrices $S_\pm(u)$, is also preserved under the regularized limit. That is 
\begin{thm}\cite{Xu}
Theorem \ref{quantumgpucat} is true for $u=u_{\rm cat}$.
\end{thm}
Similarly, we define the coproduct of the Stokes matrices at $u_{\rm cat}$ by \eqref{coS} replacing $u$ by $u_{\rm cat}$.

\section{Crystals and combinatorics of Young tableaux from the WKB approximation}

\subsection{WKB approximation of quantum Stokes matrices at caterpillar points}\label{sec:WKB}
Since the derivative in the system \eqref{texteq} is multiplied by a small parameter $1/h$, we will call the leading term, as $h\rightarrow +\infty$, of Stokes matrices of \eqref{texteq} as the WKB approximation. In this subsection, we study the WKB approximation of Stokes matrices at $u_{\rm cat}$ given in \eqref{explicitS}.

First, for any Gelfand-Tsetlin pattern $\Lambda\in GT_n(\lambda)$ given in Section \ref{GZbasis}, set
\begin{align}\label{capx}
X^{(k)}_j(\Lambda)&:=\sum_{i=1}^j(-\lambda^{(k)}_i+\lambda^{(k-1)}_{i-1}-\lambda^{(k)}_{i-1}+\lambda^{(k+1)}_{i}), \hspace{3mm} 1\le j\le k+1,\\
Y^{(k)}_j(\Lambda)&:=\sum_{i=j}^k(\lambda^{(k)}_i-\lambda^{(k-1)}_i+\lambda^{(k)}_{i+1}-\lambda^{(k+1)}_{i+1}), \hspace{3mm} 0\le j\le k,
\end{align}
where by convention $\lambda^{(k)}_j=0$ if not $1\le j\le k$. Then we define 
\begin{align}
wt_k(\Lambda)&=\sum_{i=1}^k\lambda^{(k)}_i(\Lambda)-\sum_{i=1}^{k-1}\lambda^{(k-1)}_i(\Lambda),\\
\label{maximal}
\varepsilon_k(\Lambda)&={\rm max}\{X_1^{(k)}(\Lambda), X_2^{(k)}(\Lambda),...,X_k^{(k)}(\Lambda)\},\\
\phi_k(\Lambda)&={\rm max}\{Y_1^{(k)}(\Lambda), Y_2^{(k)}(\Lambda),...,Y_k^{(k)}(\Lambda)\}.
\label{minimal}
\end{align}
Note that 
\begin{equation}\label{need}
    \phi_k(\Lambda)-\varepsilon_k(\Lambda)=wt_k(\Lambda)-wt_{k+1}(\Lambda)=Y^{(k)}_j(\Lambda)-X^{(k)}_j(\Lambda) \text{ for all } j=0,...,k+1.
\end{equation}
Furthermore, we define the functions $l_1(\Lambda)<\cdots<l_{m_k}(\Lambda)$ of $\Lambda\in GT_n(\lambda)$ be those ordered labels such that \begin{eqnarray}\label{l}
X^{(k)}_{l_j}(\Lambda)=\varepsilon_k(\Lambda).
\end{eqnarray}

Let us denote by 
\begin{equation}\label{PkGZ}
    ^k{GT}_n(\lambda):=\{\Lambda\in GT_n(\lambda)~|~\Lambda-\delta^{(k)}_{l_{m_k}}\in GT_n(\lambda) \ \text{and} \ \Lambda-\delta^{(k)}_{l_i}\notin GT_n(\lambda) \ \text{for} \ i=1,...,m_k-1\}.
\end{equation}
In particular, all $\Lambda$ satisfying ${m_k}=1$ (i.e., there is only one index $l_1$ such that $X^{(k)}_{l_1}(\Lambda)=\varepsilon_k(\Lambda)$) are in $^kGT_n(\lambda)$. The patterns in $^kGT_n(\lambda)$ are called generic, since the complements are cut out by various equalities between $\lambda^{(j)}_i$. 
Let $\{\xi_\Lambda \}$ be the Gelfand-Tsetlin basis given in Section \ref{GZbasis}. Then by a direction computation using the formula in Theorem \ref{explicitS}, we have (see \cite{Xu} for a similar computation)
\begin{pro}\label{qleadingf}
For $k=1,...,n-1$ and for any $\Lambda\in {^kGT}_n(\lambda)$, there exist real valued functions $\psi_{k}(\Lambda)$ of the patterns $\Lambda$ (independent of $q=e^{h/2}$) such that
\begin{equation*}
S_{+}(u_{\rm cat})^{L(\lambda)}_{k,k+1}\cdot \xi_\Lambda \sim q^{\phi_k(\Lambda)+\I \psi_{k}(\Lambda)} \xi_{\Lambda-\delta^{(k)}_{l_{m_k}(\Lambda)}}, \hspace{3mm} \text{as} \ \ q\rightarrow \infty.
\end{equation*}
Here recall the pattern $\Lambda-\delta^{(k)}_{l_{m_k}}$ is obtained from $\Lambda$ by replacing $\lambda^{(k)}_{l_{m_k}}$ by $\lambda^{(k)}_{l_{m_k}}- 1$. 
\end{pro}

\begin{rmk}
We remark that for the nongeneric basis element $\xi_\Lambda$, that is $m_k(\Lambda)>1$ and some $\Lambda-\delta^{(k)}_{l_{i}}\in GT_n(\lambda)$ for $i> 1$, the behaviour is 
\[S_{+}(u_{\rm cat})^{L(\lambda)}_{k,k+1}\cdot \xi_\Lambda\sim q^{\phi_k(\Lambda)+\I \beta_{k}} \xi_{\Lambda-\delta^{(k)}_{l_1(\Lambda)}}+\cdots+q^{\phi_k(\Lambda)+\I \psi_{k}} \xi_{\Lambda-\delta^{(k)}_{l_{m_k}(\Lambda)}}, \]
where the phases $\psi_k$ and $\beta_{k}$ are different, therefore the WKB leading term is not well defined. 
\end{rmk}

Similarly, let us denote by 
\begin{equation}
    GT_n^k(\lambda):=\{\Lambda\in GT_n(\lambda)~|~\Lambda+\delta^{(k)}_{l_1}\in GT_n(\lambda) \ \text{and} \ \Lambda+\delta^{(k)}_{l_i}\notin GT_n(\lambda) \ \text{for} \ i=2,...,m_k\}.
\end{equation}

\begin{pro}\label{qleading}
For $k=1,...,n-1$ and for any $\Lambda\in GT^k_n(\lambda)$, there exist real valued functions $\theta_{k}(\Lambda)$ of the patterns $\Lambda$ (independent of $q=e^{h/2}$) such that
\begin{equation}\label{qexpansion}
S_{-}(u_{\rm cat})^{L(\lambda)}_{k+1,k}\cdot \xi_\Lambda \sim q^{\varepsilon_k(\Lambda)+wt_{k+1}(\Lambda)-wt_{k}(\Lambda)+\I \theta_{k}} \xi_{\Lambda+\delta^{(k)}_{l_1(\Lambda)}} , \hspace{3mm} \text{as} \ \ q\rightarrow \infty.
\end{equation}
Here recall the pattern $\Lambda+\delta^{(k)}_{l_1}$ is obtained from $\Lambda$ by replacing $\lambda^{(k)}_{l_1}$ by $\lambda^{(k)}_{l_1}+ 1$. 
\end{pro}

Now, we can introduce the WKB operators $\widetilde{e_k}$ and $\widetilde{f_k}$, encoding respectively the leading asymptotics of the entries $S_{-}(u_{\rm cat})^{L(\lambda)}_{k+1,k}$ and $S_{+}(u_{\rm cat})^{L(\lambda)}_{k,k+1}$ of the Stokes matrices.
\begin{rmk}
Although the operators $\widetilde{e_k}$ and $\widetilde{f_k}$ are only defined for generic basis vectors, there is a unique "continuous extension" of $\widetilde{e_k}$ and $\widetilde{f_k}$ to the whole $GT_n(\lambda)$, see \cite{Xu} for more details. 
\end{rmk}
\begin{defi}\label{CWKBe}
For each $k$, the WKB operators $\widetilde{e_k}$ and $\widetilde{f_k}$ from $E_{GT}(\lambda)$ to $E_{GT}(\lambda)\cup \{0\}$ are given by
\begin{eqnarray}\label{WKBoperator}
\widetilde{e_k}\cdot \xi_{\Lambda}:=\xi_{\Lambda+\delta^{(k)}_{l_1(\Lambda)}}, \text{  and  } \widetilde{f_k}\cdot\xi_{\Lambda}:=\xi_{\Lambda-\delta^{(k)}_{l_{m_k}(\Lambda)}}.
\end{eqnarray}
Here recall that $l_1(\Lambda)$ is the integer given by \eqref{l}, and it is supposed that $\widetilde{e_k}\cdot \xi_\Lambda$ is zero if $\Lambda+\delta^{(k)}_{l_1(\Lambda)}$ doesn't belong to $P_{GT}(\lambda;{\mathbb{Z}})$.
\end{defi}

We call $(E_{GT}(\lambda),\widetilde{e_k},\widetilde{f_k},\varepsilon_k,\phi_k)$ the WKB datum of the Stokes matrices $S_{\pm}(u_{\rm cat})^{L(\lambda)}$ at the caterpillar point. 
\begin{thm}\cite{Xu}\label{SWKB}
The WKB datum $(E_{GT}(\lambda),\widetilde{e_k},\widetilde{f_k},\varepsilon_k,\phi_k)$ is a $\frak{gl}_n$-crystal, that is isomorphic to the $\frak{gl}_n$-crystal $B_n(\lambda)$ under the map $\xi_\Lambda\mapsto \tau(\Lambda)\in B_n(\lambda)$.
\end{thm}

\section{The Robinson–Schensted correspondence from the WKB approximation of the connection matrix}\label{YfromWKB}

\subsection{Deletion operator from the WKB approximation of the connection matrix}\label{SfromWKB}

Let us take a partition $\lambda=(\lambda_1,...,\lambda_n)$, seen as a Young diagram.
The decomposition of the tensor product of the representation $L(\lambda)$ and the natural representation $V=\mathbb{C}^n$ is
\begin{equation}
    L(\lambda)\otimes \mathbb{C}^n\cong \underset{\mu}{\bigoplus} L(\mu)
\end{equation}
where $\mu$ runs over Young diagrams obtained by adding one node to $\lambda$. 

Let us assume $\lambda\stackrel{a}{\longleftarrow} \mu$, i.e., the node is added at the $a$-th row. 
Then the embedding $L(\mu)\subset  L(\lambda)\otimes \mathbb{C}^n$ is described by Wigner coefficients: given a basis vector $\xi_{\Lambda}\in L(\mu)$ with $\Lambda\in GT_n(\mu)$,
for $(i_n,...,i_j)$ such that $i_n=a$ and $1\le i_k\le k$ $(j\le k\le n)$, we define a GT pattern $(\Lambda; i_n,...,i_j)\in GT_n(\lambda)$ by 

\begin{equation}
   \lambda^{(k)}_i(\Lambda; i_n,...,i_j) =\left\{
          \begin{array}{lr}
            \lambda^{(k)}_{i_k}(\Lambda)-1,   & \text{if} \ \ j\le k\le n \ \text{ and } \ i=i_k, \\
            \lambda^{(k)}_{i_k}(\Lambda), & \text{otherwise}.
             \end{array}
\right.
\end{equation}
Then the embedding $\xi_{\Lambda}\in L(\mu)\subset  L(\lambda)\otimes \mathbb{C}^n$ reads as
\begin{equation}
\xi_{\Lambda}=\sum_{j=1}^n\sum_{i_n=a, i_{n-1},...,i_j} \omega(\Lambda; i_n,...,i_j)\times \xi_{(\Lambda; i_n,...,i_j) } \otimes v_j,   
\end{equation}
where the Wigner coefficients are
\begin{equation}
    \omega(\Lambda; i_n,...,i_j)=\omega^{(1)}\left(\begin{array}{c}
\lambda^{(j)}_1, ..., \lambda^{(j)}_j \\
\lambda^{(j-1)}_1, ..., \lambda^{(j-1)}_{j-1}
\end{array} ; \begin{array}{c}
i_j \\

\end{array} \right)\prod_{k=j+1}^n \omega^{(2)}\left(\begin{array}{c}
\lambda^{(k)}_1, ..., \lambda^{(k)}_k \\
\lambda^{(k-1)}_1, ..., \lambda^{(k-1)}_{k-1}
\end{array} ; \begin{array}{c}
i_k \\
i_{k-1}
\end{array} \right)
\end{equation}
Here for a sequence of numbers $f$
\begin{align}
&\omega^{(1)}\left(\begin{array}{c}
\lambda^{(j)}_1, ..., \lambda^{(j)}_j \\
\lambda^{(j-1)}_1, ..., \lambda^{(j-1)}_{j-1}
\end{array} ; \begin{array}{c}
i_j \\
\end{array} \right)=\sqrt{\frac{\prod_{t\le j-1}(l^{(j-1)}_{t}-l_{i_j}^{(j)})}{\prod_{t\le j, t\ne i_j}(l^{(j)}_{t}-l^{(j)}_{i_j}+1)}}\\
&\omega^{(2)}\left(\begin{array}{c}
\lambda^{(k)}_1, ..., \lambda^{(k)}_k \\
\lambda^{(k-1)}_1, ..., \lambda^{(k-1)}_{k-1}
\end{array} ; \begin{array}{c}
i_k \\
i_{k-1}
\end{array} \right)=S(i_{k-1}-i_k) \sqrt{\prod_{t\le k, t\ne i_k}\frac{l^{(k)}_{t}-l^{(k-1)}_{i_{k-1}}+1}{l^{(k)}_{t}-l^{(k)}_{i_{k}}+1} \prod_{t\le k-1, t\ne i_{k-1}}\frac{l^{(k-1)}_{t}-l^{(k)}_{i_k}}{l^{(k-1)}_{t}-l^{(k-1)}_{i_{k-1}}}}.
\end{align}
Here 
\[l^{(i)}_{t}=\lambda^{(i)}_{t}(\Lambda)-t \ \ \text{  for any } 1\le t\le i\le n,\]
and 
\[ S(i_{k-1}-i_k)=\left\{
          \begin{array}{lr}
            1,   & \text{if} \ \ i_{k-1}-i_k\ge 0, \\
            -1, & \text{if} \ \ i_{k-1}-i_k< 0.
             \end{array}
\right. \]

\begin{thm}\label{WKBC}
Given a basis vector $\xi_{\Lambda}\in L(\mu)\subset  L(\lambda)\otimes \mathbb{C}^n$ with $\lambda\stackrel{a}{\longleftarrow} \mu$,
there exists real valued function $\theta(\Lambda)$ such that 
\begin{eqnarray*}
\mathop{\rm lim}\limits_{h\rightarrow +\infty} C(u_{\rm cat})^{L(\lambda)}\cdot e^{\I h \theta(\Lambda)} \xi_\Lambda =\xi_{\leftarrow_a\Lambda}\otimes v_j.
\end{eqnarray*}
Here $\xi_{\leftarrow_a\Lambda}:=\tau^{-1}(\leftarrow_a\tau(\Lambda))\in L(\lambda)$ (the $\tau^{-1}:B_n(\lambda) \rightarrow GT_n(\lambda)$ is given in Definition \ref{tau}), and $\leftarrow_a\tau(\Lambda)\in B_n(\lambda)$ is the semistandard tableaux obtained by the deletion operation as in Definition \ref{delop} on the semistandard tableaux $\tau(\Lambda)\in B_n(\mu)$ (with the map $\tau:GT_n(\mu)\rightarrow B_n(\mu)$ with respect to the partition $\mu$).
\end{thm}

\begin{proof}

Given the $\xi_{\Lambda}\in L(\mu)\subset  L(\lambda)\otimes \mathbb{C}^n$ with $\lambda\stackrel{a}{\longleftarrow} \mu$, let us write 
\[C(u_{\rm cat})^{L(\lambda)}=\sum_{i,j=1,...,n}C(u_{\rm cat})^{L(\lambda)}_{ij}\otimes E_{ij}\in {\rm End}(L(\lambda))\otimes {\rm End}(\mathbb{C}^n),\] then via the Wigner coefficients
\begin{equation}\label{Caction}
    C(u_{\rm cat})^{L(\lambda)}\cdot \xi_\Lambda=\sum_{i=1}^n\sum_{j=1}^n\sum_{i_{n-1},...,i_j} \omega(m;i_n,...,i_j)\times \Big(\overrightarrow{\prod_{k=2,...,n}}\widetilde{C^{(k)}}\Big)_{ik}\cdot P_{kj}\cdot \xi_{(\Lambda;i_n,...,i_j)})\otimes v_i.
\end{equation}
Let us first use the following lemma to simply \eqref{Caction}.
\begin{lem}\label{minorwigner}
We have (recall that the index $i_n=a$)
\begin{equation}\label{Qwig1}
\sum_{j=1}^n  \sum_{i_{n-1},...,i_j}  \omega(\Lambda; i_n,...,i_j)\times  P_{aj}\cdot \xi_{(\Lambda; i_n,...,i_j) } =\xi_{(\Lambda; i_n=a) },
\end{equation}
and for all the other indices $k\ne a$ and $j$
\begin{equation}\label{Qwig2}
\sum_{j=1}^n\sum_{i_n=a, i_{n-1},...,i_j}\omega(\Lambda; i_n,...,i_j)\times  P_{kj}\cdot \xi_{(\Lambda; i_n,...,i_j) } =0,
\end{equation}
\end{lem}
\begin{proof}
Recall that the entries of the matrix $Q$ are quantum minors of $\mathcal{T}$. The action of these quantum minors on the Gelfand-Tsetlin basis are given by Wigner coefficients. In particular,
\begin{equation}\label{PQ1}
    Q_{j i_n}\cdot \xi_{(\Lambda; i_n) } =  \sum_{i_{n-1},...,i_j} \omega(\Lambda; i_n,...,i_j)\times   \xi_{(\Lambda; i_n,...,i_j) }.
\end{equation}
By Lemma \ref{diagT}, we have $P_{kj}Q_{ji_n}=\delta_{k,i_n}$. Together with \eqref{PQ1}, it proves the identities \eqref{Qwig1} and \eqref{Qwig2}.
\end{proof}

Recall that the $(i,i)$-entry of $\widetilde{C^{(k)}}$ is the identity in ${\rm End}(L(\lambda))$ for all $1\le k\le i-1$. It follows from \eqref{Caction} and Lemma \ref{minorwigner} that 
\begin{align}\nonumber 
C(u_{\rm cat})^{L(\lambda)}_{ik}\cdot \xi_{\Lambda}&=\sum_{i=1}^n\left(\Big(\overrightarrow{\prod_{k=2,...,n}}\widetilde{C^{(k)}}\Big)_{ia}\cdot\xi_{(\Lambda; i_n=a) }\right)\otimes v_i\\ \label{expansion}
&=\sum_{i=1}^n\left(\sum_{t_{i},...,t_{n-1}} \widetilde{C^{(i)}_{it_i}}\widetilde{C^{(i+1)}_{t_it_{i+1}}}\cdots \widetilde{C^{(n)}_{t_{n-1}, a}}\cdot \xi_{(\Lambda; i_n=a) }\right)\otimes v_i,
\end{align}
where the summation is over all possible $t_l=1,...,l$ with $l=i,...,n-1$. 

In order to determine the leading asymptotics on the right hand side of \eqref{expansion} as $h\rightarrow+\infty$, 
we need the following lemma.

\begin{lem}\label{wkbrelC}
For any basis vector $\xi_{\Lambda'}$ with $\Lambda'\in GT_n(\lambda)$, we have as $h\rightarrow+\infty$ (with $q=e^{\frac{h}{2}}$)
\begin{align}
 \widetilde{{C}^{(k+1)}_{ij}}\cdot \xi_{\Lambda'}&\sim q^{i-j+1+\sum_{i<l<j}(\lambda^{(k)}_{l}-\lambda^{(k+1)}_l)}q^{\I \theta_{ij}}\times \xi_{\Lambda'-\delta^{(k)}_{i}}, & \text{if } \ i<j,\\
       \widetilde{{C}^{(k+1)}_{ij}}\cdot \xi_{\Lambda'}&\sim q^{i-j+\lambda_i^{(k)}-\lambda_j^{(k+1)}+\sum_{l=j}^{i-1}({\lambda_{l+1}^{(k+1)}}-{\lambda_l^{(k)}})}q^{\I \theta_{ij}}\times \xi_{\Lambda'-\delta^{(k)}_{i}}, &\text{if } \ i\ge j \text{ and } i\ne n,\\
       \widetilde{{C}^{(k+1)}_{ij}}\cdot \xi_{\Lambda'}&\sim q^{\sum_{1\le l<j}{(\lambda^{(k)}_l}-\lambda^{(k+1)}_l)}q^{\I \theta_{ij}}\times \xi_{\Lambda'}, & \text{if } \ i=n,
\end{align}
where $\theta_{ij}$ is certain real function depending on the pattern $\Lambda'$, and $\lambda^{(i)}_j=\lambda^{(i)}_j({\Lambda'}).$
    
\end{lem}
\begin{proof}
The asymptotics of gamma function states that
for $r$ a real number, as $h\rightarrow+\infty$
\begin{equation}\label{gamma}
{\rm In}\left(\Gamma\left(1+\frac{rh}{2\pi \I }\right)\right)\sim \frac{rh{\rm In}(h)}{2\pi \I }+\frac{rh}{2\pi \I }{\rm In}\left(\frac{|r|}{2\pi}\right)- \frac{|r|h}{4}-\frac{rh}{2\pi \I }+\frac{1}{2}{\rm In}\left(\frac{rh}{\I }\right). 
\end{equation}
By \eqref{eigenvalues}, \[\zeta^{(k)}_i\cdot \xi_{\Lambda'} =(\lambda^{(k)}_i({\Lambda'})-i+1)\cdot\xi_{\Lambda'}.\] 

The lemma follows directly from \eqref{gamma}, the interlacing inequalities between $\lambda^{(i)}_j$ for $i=k-1,k,k+1$ and the expressions in Proposition \ref{GZaction}. \end{proof}

\begin{cor}\label{leadC}
As $h\rightarrow+\infty$, the real part $f_{ij}(\Lambda')$ in the exponent of the leading term (with $q=e^{\frac{\hbar}{2}}$)
\[\widetilde{C^{(k+1)}_{ij}}\cdot \xi_{\Lambda'}\sim q^{f_{ij}(\Lambda')+\I \theta_{ij}(\Lambda')}\times \xi_{\Lambda''}\]
satisfies that 
\begin{equation}\label{real0}
    f_{ij}(\Lambda')=\left\{
          \begin{array}{lr}
            0,   & \text{if } \ i=j\le k \ \text{ and } \lambda^{(k+1)}_i(\Lambda')=\lambda^{(k)}_i(\Lambda') \\
           0, & \text{if } \ i=j-1  \\
           0, & \text{if } \ i=1 \text{ and } j=k+1 \\
           <0, & \text{otherwise} 
             \end{array}\right.
\end{equation}
\end{cor}
Using Corollary \eqref{leadC}, we check that in the summation \eqref{expansion},
there exists a unique sequence of indices $i_n=a, i_{n-1},...,i_j$ such that the summand 
\[\widetilde{C^{(j+1)}_{i_ji_{j+1}}} \widetilde{C^{(j+2)}_{i_{j+1}i_{j+2}}} \cdots \widetilde{C^{(n)}_{i_{n-1}i_n}}\cdot \xi_{(\Lambda;i_n=a)}\]
is the only term whose real part in the exponent of the asymptotics as $q=e^{\frac{h}{2}}\rightarrow+\infty$
is $0$. Furthermore, using Lemma \ref{wkbrelC} we have for the unique sequence $i_n=a, i_{n-1},...,i_j$
\[\lim_{h\rightarrow+\infty}\widetilde{C^{(j+1)}_{i_ji_{j+1}}} \widetilde{C^{(j+2)}_{i_{j+1}i_{j+2}}} \cdots \widetilde{C^{(n)}_{i_{n-1}i_n}}\cdot q^{\I \theta(\Lambda;i_n)} \xi_{(\Lambda;i_n=a)}=\xi_{(\Lambda;i_n,...,i_j)}),\]
for some real number $\theta(\Lambda;i_n)$. In fact, the unique sequence $i_n=a,i_{n-1},...,i_j$ can be defined by induction: assume that $i_n,...,i_k$ are defined, then set $\Lambda'=(\Lambda;i_n,...,i_k)$ in Corollary \ref{leadC}, by \eqref{real0} we simply proceed as follows
\begin{itemize}
    \item[(1)] If $k=1$, then set $j=1$ and stop.

    \item[(2)] If $k>1$ and $\lambda^{(k)}_{i_k}=\lambda^{(k)}_{i_{k-1}}$, then set $i_{k-1}=i_k$.

    \item[(3)] If $k>1$, $\lambda^{(k)}_{i_k}>\lambda^{(k)}_{i_{k-1}}$ and $i_k>1$, then set $i_{k-1}=i_k-1$.

    \item[(4)] If $k>1$, $\lambda^{(k)}_{i_k}>\lambda^{(k)}_{i_{k-1}}$ and $i_k=1$, then set $j=k$ and stop.
\end{itemize}

We check that the above procedure $(1)-(4)$ is equivalent to the deletion operation $\leftarrow_a\Lambda$ in Definition \ref{delop}. That is the unique sequence $i_n,...,i_j$, given by procedure $(1)-(4)$ with initial $i_n=a$, is such that $(\Lambda;i_n,...,i_j)=\leftarrow_a\Lambda$. It thus finishes the proof.
\end{proof}

\subsection{The Robinson–Schensted correspondence from the WKB approximation}
By repeatedly utilizing Theorem \ref{WKBC}, we get
    
\begin{thm}
If $\xi_\Lambda \in V^{\otimes N}=\underset{SY}{\bigoplus} L(SY_R)$ is a vector parametrized by the pair of Young tableau $RS(\omega)=(P(\omega), Q(\omega))$ with the word $\omega=w_1\cdots w_N$, then there exits a real number $\theta(\Lambda)$ such that 
\begin{eqnarray*}
\mathop{\rm lim}\limits_{h\rightarrow +\infty}\prod_{k=1}^{N-1} \left({\rm Id}_{N-k-1}\otimes C(u_{\rm cat})^{V^{\otimes (N-k)}}\right)\cdot q^{\I \theta(\Lambda)}\xi_\Lambda=v_{w_1}\otimes \cdots \otimes v_{w_N}.
\end{eqnarray*}
\end{thm}

\subsection{Tensor products of crystals from the WKB approximation of coproduct}\label{sec:tensor} 
Now given two representations $L(\mu)$ and $L(\nu)$, we consider the actions of the operators $\Delta(S_{+}(u_{\rm cat})_{k,k+1})$ and $\Delta(S_{-}(u_{\rm cat})_{k+1,k})$ (as defined in \eqref{coS}) on the tensor product $L(\mu)\otimes L(\nu)$. Let $\{\xi_{\Lambda_1}\}$ and $\{\xi_{\Lambda_2}\}$ be the basis of $L(\lambda_1)$ and $L(\lambda_2)$ respectively. Similar to Proposition \ref{qleading}, one can compute the WKB leading term of the operators defined in \eqref{copro1}-\eqref{copro2} under the basis $\xi_{\Lambda_1}\otimes \xi_{\Lambda_2}$. 
\begin{pro}\label{coproduct}
The WKB approximation of the operators $\Delta(S_{+}(u_{\rm cat})_{k,k+1})$ and $\Delta(S_{-}(u_{\rm cat})_{k+1,k})$ induce the crystal operators $\widetilde{e_k}$ and $\widetilde{f_k}$ on the tensor product $E_{GT}(\mu)\otimes E_{GT}(\nu)$ of $\frak{gl}_n-$crystals respectively. 
\end{pro}

\begin{proof} Take any vector $\xi_{\Lambda_1}\otimes \xi_{\Lambda_2}\in L(\mu)\otimes L(\nu)$. Assume that $\xi_{\Lambda_1}$ and $\xi_{\Lambda_2}$ are generic as in the sense of Propositions \ref{qleading} and \ref{qleadingf}. As $q=e^{\frac{h}{2}}\rightarrow +\infty$,
\begin{align*}
&\Delta(S_{+}(u_{\rm cat})_{k,k+1})^{L(\mu)\otimes L(\nu)}\cdot (\xi_{\Lambda_1}\otimes \xi_{\Lambda_2})\\
=&\Big(q^{e_{ii}-e_{i+1,i+1}}\cdot \xi_{\Lambda_1}\Big)\otimes \Big( S_{+}(u_{\rm cat})^{L(\nu)}_{k,k+1}\cdot \xi_{\Lambda_2}\Big)+\Big( S_{+}(u_{\rm cat})_{k,k+1}\cdot \xi_{\Lambda_1}\Big)\otimes \Big(1\cdot \xi_{\Lambda_2}\Big)\\
\sim & q^{wt_k(\Lambda_1)-wt_{k+1}(\Lambda_1)+\phi_k(\Lambda_2)+\I \theta_{k}(\Lambda_2)} \xi_{\Lambda_1}\otimes \xi_{\Lambda_2-\delta^{(k)}_{l_{m_k}(\Lambda_2)}}+q^{\phi_k(\Lambda_1)+\I \theta_{k}(\Lambda_1)} \xi_{\Lambda_1-\delta^{(k)}_{l_{m_k}(\Lambda_1)}}\otimes \xi_{\Lambda_2}.
\end{align*}

To determine which term is the leading one, we only need to compare the exponents. Using the identity \eqref{need}, it gives
\begin{equation}
\mathcal{WKB}\left( \Delta(S_{+}(u_{\rm cat})_{k,k+1})^{L(\mu)\otimes L(\nu)}\cdot (\xi_{\Lambda_1}\otimes \xi_{\Lambda_2})\right)=\left\{
          \begin{array}{lr}
            \xi_{\Lambda_1}\otimes \xi_{\Lambda_2-\delta^{(k)}_{l_{m_k}(\Lambda_2)}},   & \text{if } \phi_k(\Lambda_2)>\varepsilon_k(\Lambda_1)  \\
           \xi_{\Lambda_1-\delta^{(k)}_{l_{m_k}(\Lambda_1)}}\otimes \xi_{\Lambda_2}, & \text{otherwise}
             \end{array}\right.
\end{equation}
It is compatible with the definition of the operator $\widetilde{f_k}$ on the tensor product given in Definition \ref{crytensor}. 

Similarly, as $q=e^{\frac{h}{2}}\rightarrow +\infty$,
\begin{align*}
&\Delta(S_{-}(u_{\rm cat})_{k+1,k})^{L(\mu)\otimes L(\nu)}\cdot (\xi_{\Lambda_1}\otimes \xi_{\Lambda_2})\\
=&\Big(q^{e_{i+1,i+1}-e_{ii}}\cdot \xi_{\Lambda_1}\Big)\otimes \Big( S_{h-}(u_{\rm cat})^{L(\nu)}_{k+1,k}\cdot \xi_{\Lambda_2}\Big)+\Big( S_{h-}(u_{\rm cat})_{k+1,k}\cdot \xi_{\Lambda_1}\Big)\otimes \Big(q^{2e_{i+1,i+1}-2e_{ii}}\cdot \xi_{\Lambda_2}\Big)\\
\sim & q^{-wt_{k}(\Lambda_1)+wt_{k+1}(\Lambda_1)+\varepsilon_k(\Lambda_2)+wt_{k+1}(\Lambda_2)-wt_{k}(\Lambda_2)+\I \theta_{k}(\Lambda_2)} \xi_{\Lambda_1}\otimes \xi_{\Lambda_2+\delta^{(k)}_{l_1(\Lambda_2)}}\\
&+q^{\varepsilon_k(\Lambda_1)+wt_{k+1}(\Lambda_1)-wt_{k}(\Lambda_1)+2wt_{k+1}(\Lambda_2)-2wt_{k}(\Lambda_2)+\I \theta_{k}(\Lambda_1)} \xi_{\Lambda_1+\delta^{(k)}_{l_1(\Lambda_1)}}\otimes \xi_{\Lambda_2}.
\end{align*}
Comparing the exponents gives (using \eqref{need})
\begin{equation}
\mathcal{WKB}\left( \Delta(S_{-}(u_{\rm cat})_{k+1,k})^{L(\mu)\otimes L(\nu)}\cdot (\xi_{\Lambda_1}\otimes \xi_{\Lambda_2})\right)=\left\{
          \begin{array}{lr}
            \xi_{\Lambda_1}\otimes \xi_{\Lambda_2+\delta^{(k)}_{l_1(\Lambda_2)}},   & \text{if } \phi_k(\Lambda_2)>\varepsilon_k(\Lambda_1)  \\
           \xi_{\Lambda_1+\delta^{(k)}_{l_1(\Lambda_1)}}\otimes \xi_{\Lambda_2}, & \text{otherwise}
             \end{array}\right.
\end{equation}
It shows that the WKB operator $\mathcal{WKB}\left( \Delta(S_-)_{k+1,k}\cdot (\xi_{\Lambda_1}\otimes \xi_{\Lambda_2})\right)$ is compatible with the definition of the crystal operator $\widetilde{e_k}$ on the tensor product given in Definition \ref{crytensor}. 
\end{proof}

\subsection{Littlewood–Richardson rule and the WKB approximation of coproduct}\label{LRWKB}

There is an irreducible decomposition of the tensor product of two irreducible representations $L(\mu)$ and $L(\nu)$ into $\bigoplus_\lambda m^\lambda_{\mu\nu} L(\lambda)$. Following Theorem \ref{LRc}, we introduce the map 
\begin{align}\label{mapLR}
LR:~L(\mu)\otimes L(\nu) &\rightarrow \bigoplus_\lambda L(\lambda)\times LR^\lambda_{\mu\nu} \\
\xi_{\Lambda_1}\otimes \xi_{\Lambda_2} &\mapsto \ \ (\xi_{\tau^{-1}(P(\tau(\Lambda_1), \tau(\Lambda_2)))}, Q(\tau(\Lambda_1), \tau(\Lambda_2))),
\end{align}
where $\xi_{\Lambda_1}$ and $\xi_{\Lambda_2}$ are the GT basis vectors of $L(\mu)$ and $L(\nu)$ with respect to the GT patterns $\Lambda_1\in GT_n(\mu)$ and $\Lambda_2\in GT_n(\nu)$ respectively, and then $\tau(\Lambda_1)\in B_n(\mu)$, $\tau(\Lambda_2)\in B_n(\nu)$ as in Definition \ref{tau}.

\begin{thm}\label{wkbcoproduct}
Let $\{\xi_{\Lambda_1}\}$ and $\{\xi_{\Lambda_2}\}$ be any generic basis of $L(\mu)$ and $L(\nu)$, parameterized by $\Lambda_1\in GT_n(\mu)$ and $\Lambda_2\in GT_n(\nu)$ respectively, then
\begin{eqnarray*}
LR\left(\mathcal{WKB}\left( \Delta(S_+(u_{\rm cat}))^{L(\mu)\otimes L(\nu)}_{k,k+1}\cdot (\xi_{\Lambda_1}\otimes \xi_{\Lambda_2})\right)\right)=\mathcal{WKB}\left( (S_+(u_{\rm cat})^{L(\mu)\otimes L(\nu)})_{k,k+1}\cdot LR(\xi_{\Lambda_1}\otimes \xi_{\Lambda_2}) \right),
\end{eqnarray*}
where the map $LR:L(\mu)\otimes L(\nu)\rightarrow \underset{\lambda}{\bigoplus} L(\lambda)\times LR^\lambda_{\mu\nu}$ is given by \eqref{mapLR} according to the Littlewood–Richardson rule.
\end{thm}
\begin{proof}
Firstly, the action of $S_+(u_{\rm cat})_{k,k+1}^{L(\mu)\otimes L(\nu)}\in {\rm End}(L(\mu)\otimes L(\nu))$ can be restricted to each component $L(\lambda)$. By Theorem \ref{SWKB}, the WKB approximation of $S_+(u_{\rm cat})_{k,k+1}^{L(\mu)\otimes L(\nu)}$ on the component $L(\lambda)$ gives rise to a $\frak{gl}_n$-crystal $E_{GT}(\lambda)$ that is isomorphic to $B_n(\lambda)$. Therefore, we have
\begin{equation}\label{WKBLR}
     \mathcal{WKB}\left( (S_+(u_{\rm cat})^{L(\mu)\otimes L(\nu)})_{k,k+1}\cdot LR(\xi_{\Lambda_1}\otimes \xi_{\Lambda_2}) \right)=\widetilde{f_k}\left(LR(\xi_{\Lambda_1}\otimes \xi_{\Lambda_2})\right).
\end{equation}
Here $\widetilde{f_k}$ is the crystal operator on $E_{GT}(\lambda)\cong B_n(\lambda)$.

Secondly, by Proposition \ref{coproduct},
\begin{equation}
    \mathcal{WKB}\left( \Delta(S(u_{\rm cat})_+)^{L(\mu)\otimes L(\nu)}_{k,k+1}\cdot (\xi_{\Lambda_1}\otimes \xi_{\Lambda_2})\right)= \widetilde{f_k} (\xi_{\Lambda_1}\otimes \xi_{\Lambda_2}).
\end{equation}
Here $\widetilde{f_k}$ is the crystal operator on the tensor $E_{GT}(\mu)\otimes E_{GT}(\nu)$ of crystals.

Thirdly, by Theorem \ref{LRcrystal}, we have
\begin{equation}\label{LRLR}
\widetilde{f_k}\left(LR(\xi_{\Lambda_1}\otimes \xi_{\Lambda_2})\right)=LR\left(\widetilde{f_k}\left(\xi_{\Lambda_1}\otimes \xi_{\Lambda_2}\right)\right).
\end{equation}
Combining the above three identities \eqref{WKBLR}-\eqref{LRLR}, we get the desired identity.
\end{proof}

\Addresses
\end{document}